\documentclass[11pt,reqno]{amsart}

\usepackage{marginnote,amssymb,mathrsfs,latexsym,verbatim,pdfsync,color,mathabx}

\newcommand{\bbC}{{\mathbb C}}

\newcommand{\bbH}{{\mathbb H}}

\newcommand{\bbN}{{\mathbb N}}
\newcommand{\bbR}{{\mathbb R}}

\def\bbC{{\mathbb C}}

\def\bbH{{\mathbb H}}

\def\bbN{{\mathbb N}}

\def\bbR{{\mathbb R}}

\def\bX{{\mathbf X}}

\def\cN{{\mathcal N}}

\def\cS{{\mathcal S}}

\def\cX{{\mathcal X}}

\def\sA{{\mathscr A}}
\def\sB{{\mathscr B}}

\def\sF{{\mathscr F}}

\def\sI{{\mathscr I}}
\def\sJ{{\mathscr J}}

\def\sS{{\mathscr S}}

\def\eps{\varepsilon}
\def\z{\zeta} 
\def\l{\lambda} 
\def\vp{\varphi}

\def\p{\partial}

\def\ms{\medskip}

\def\endpf{\medskip\hfill $\Box$

\ms

}  
\def\epf{\endpf}

\DeclareMathOperator{\supp}{supp}
\DeclareMathOperator{\oD}{D}

\def\Dc{\Delta_\chi}

\DeclareFontFamily{U}{mathx}{\hyphenchar\font45}
\DeclareFontShape{U}{mathx}{m}{n}{
      <5> <6> <7> <8> <9> <10>
      <10.95> <12> <14.4> <17.28> <20.74> <24.88>
      mathx10
      }{}
\DeclareSymbolFont{mathx}{U}{mathx}{m}{n}
\DeclareFontSubstitution{U}{mathx}{m}{n}
\DeclareMathAccent{\widecheck}{0}{mathx}{"71}
\DeclareMathAccent{\wideparen}{0}{mathx}{"75}

\makeatletter
\newcommand{\leqnomode}{\tagsleft@true}
\newcommand{\reqnomode}{\tagsleft@false}
\makeatother

\def\itb{\item[{\tiny $\bullet$}]}

\newtheorem{thm}{Theorem}[section]

\newtheorem{lem}[thm]{Lemma}
\newtheorem{defn}[thm]{Definition}

\begin{document}

\title[Potential spaces on Lie groups]{Potential spaces on Lie groups} 

\author[T.\ Bruno]{Tommaso Bruno}
\address{Dipartimento di Scienze Matematiche ``Giuseppe Luigi Lagrange'',
  Politecnico di Torino, Corso Duca degli Abruzzi 24, 10129 Torino,
  Italy - Dipartimento di Eccellenza 2018-2022}
\email{tommaso.bruno@polito.it}

\author[M.\ M.\ Peloso]{Marco M.\ Peloso}
\address{Dipartimento di Matematica, 
Universit\`a degli Studi di Milano, 
Via C.\ Saldini 50,  
20133 Milano, Italy}
\email{marco.peloso@unimi.it}

\author[M.\ Vallarino]{Maria Vallarino}
\address{Dipartimento di Scienze Matematiche ``Giuseppe Luigi Lagrange'',
  Politecnico di Torino, Corso Duca degli Abruzzi 24, 10129 Torino,
  Italy - Dipartimento di Eccellenza 2018-2022}
\email{maria.vallarino@polito.it}

\keywords{Lie groups, Besov spaces, Triebel--Lizorkin spaces}
\thanks{{\em Math Subject Classification} 46E35, 22E30, 43A15}
\thanks{All authors are  partially supported by the grant PRIN 2015
  {\em Real and Complex Manifolds: Geometry, Topology and Harmonic
    Analysis}, and are members of the Gruppo Nazionale per l'Analisi
  Matematica, la Probabilit\`a e le loro Applicazioni (GNAMPA) of the
  Istituto Nazionale di Alta Matematica (INdAM)} 

\begin{abstract}
In this paper we discuss function spaces on a general noncompact Lie group, namely
the scales of Triebel--Lizorkin and Besov spaces, defined in terms of
a sub-Laplacian with drift.  The sub-Laplacian is written as negative
the sum of squares of a collection of left-invariant vector fields
satisfying H\"ormander's condition.  These spaces were recently
introduced by the authors.  In this paper we prove a norm characterization
in terms of finite differences, the density of test functions, and
related 
isomorphism properties.
\end{abstract}
\maketitle

\begin{centerline}
{\it Dedicated to Fulvio Ricci on the occasion of his 70th birthday} \ms
\end{centerline}

\section*{Introduction}
\ms

The theory of function spaces, regularity of integral operators, and
of solutions of differential equations, began
in the setting of Euclidean spaces, with smoothness 
measured in terms of
Sobolev and Lipschitz norms, see e.g.~\cite{Stein-SingularIntegrals}. 
A.\ Calder\'on and A.\ Zygmund developed the theory of singular
integrals, proving their  boundedness in the Lebesgue spaces, 
as well as regularity of solutions of classical
differential equations, such as the Dirichlet and Neumann problems,
in the case of a half-space and of smooth domains.  
Among the operators studied were the  singular integrals, hence
in particular the Hilbert and Riesz transforms, the Poisson integral,
and 
the heat propagator. It is worth noticing that the singularities of the integral kernels of such operators, or better, of the level sets of their moduli, were naturally
described in terms of the underlying Euclidean geometry. Such theory then included embedding and
interpolation results for Lebesgue, Sobolev and Lipschitz spaces, see
e.g.~\cite{Bergh-Lofstrom}. In this analysis, the Fourier series and transform played a crucial
role, and a noticeable application of such techniques was the
decomposition initially introduced by Littlewood and Paley, and later
developed in depth by E.\ M.\ Stein~\cite{Stein}.  The Littlewood--Paley
decomposition was initially intended to provide a substitute for the
Plancherel formula to the $L^p$-norms, with $p\neq2$, but proved to be
an invaluable tool in many other situtations.  The function spaces
that naturally arose in studying the regularity properties of
aforementioned operators, were indeed, besides the Lebesgue spaces, 
the Sobolev and Lipschitz spaces, and also the Besov spaces.  It
became then natural to obtain other characterizations for such norms,
and in this setting the Littlewood--Paley decomposition proved to be
very useful, and was also
used to define another, related, scale of spaces, the so-called
Triebel--Lizorkin spaces, see e.g.~\cite{Triebel}, which include the
Sobolev spaces as a special case.

While such theory was in its full development, L.\ H\"ormander 
produced two breakthrough results,~\cite{Hormander-translation}
and~\cite{Hormander-sum-of-squares}.   In~\cite{Hormander-translation} H\"ormander extended a previous result by Mihlin,
developing  the theory of $L^p$-multipliers of
the Laplacian. This approach also stimulated the study of a 
 class of operators that naturally appear while solving partial
 differential equations involving the Laplacian --- for instance the
 wave equation in
the Euclidean space $\bbR^d$.

 In~\cite{Hormander-sum-of-squares} H\"ormander
showed that operators that are sum
of squares of vector fields whose commutators up to a finite order
span all directions of $\bbR^d$, although non-elliptic, enjoy many
interesting and strong properties of elliptic operators, in particular
hypoelliticity.  Such phenomenon appeared for instance in the case
of the Kohn-Laplacian on the boundary of the Siegel upper half-space
in $\bbC^{d+1}$,
in the works of A.\ Kor\'anyi and S.\ V\'agi~\cite{Koranyi-Vagi},
J.\ J.\ Kohn~\cite{Kohn} and,
with most relevance to this discussion and  the present work,  of
G.\ B.\ Folland and Stein~\cite{Folland-Stein}.  The operators that
were considered in \cite{Folland-Stein}, that is the Kohn-Laplacian, the sub-Laplacian, the
so-called Folland--Stein operators, their fundamental solution, or the
relative fundamental solution in some cases, had the singularity that
could be described in terms of a different underlying geometry.
The boundary of the
Siegel upper half-space can be identified with the Heisenberg group,
and such geometry was more efficiently described using the nilpotent
Lie group structure of 
the Heisenberg group.  As metric space, the Heisenberg group $\bbH_d$
is not
equivalent to the Euclidean space $\bbR^{2d+1}$, and in fact the distance
coincides with the Carnot--Carath\'eodory distance defined by the
sub-Laplacian on  $\bbH_d$.  
The Lie algebra of $\bbH_d$ can be written as the linear span of a
family of vector fields $\bX =\{X_1,\dots,X_{2d}\}$ and of their
commutators, which reduce in fact to a single ``transversal'' vector field
$T$.  The sub-Laplacian on $\bbH_d$ is the (negative) sum of squares
$-\sum_{j=1}^{2d} X_j^2$, and thus is of the type studied by
H\"ormander in~\cite{Hormander-sum-of-squares}.
The function spaces that better describe the
smoothness of functions in this setting 
can be described by their behaviour with respect to the action of only
the vector fields $\bX$.  Such system of vector fields were called
{\em horizontal}  and they
were studied in~\cite{Folland-Stein} and~\cite{Folland-Arkiv} and
again differed from their Euclidean 
analogue.  In these papers, the authors proved
analogue of embedding and interpolation results for the newly defined 
 Sobolev and Lipschitz spaces, in the case of $\bbH_d$, and of
 Carnot--Carath\'eodory groups, respectively.\footnote{In~\cite{Folland-Arkiv}
   the Carnot--Carath\'eodory  
groups were called {\em stratified nilpotent Lie groups}.}  

These results gave tremendous impetus to the 
development of analysis on $\bbH_d$, and more in general on Carnot--Carath\'eodory 
groups.
In a series of papers, F.\ Ricci and E.\ M.\ Stein~\cite{Ricci-Stein1,
  Ricci-Stein2, Ricci-Stein3}
studied the boundedness of
singular integrals on nilpotent Lie groups, exploring again the
connection between the geometry of the metric balls, the size
properties of the integral kernels, and the boundedness of the
singular integral operators.  In~\cite{Strichartz-Hn} 
R.\ Strichartz pointed out the importance of the role of the joint
spectrum of the sub-Laplacian and $T$. In two fundamental papers, 
\cite{MRS1, MRS2} 
  D.\ M\"uller, F.\ Ricci, and E.\ M.\ Stein then proved 
the boundedness of joint spectral multipliers
of the sub-Laplacian and $T$ on $\bbH_d$ and the closely related 
Heisenberg type groups  
-- results that are yet to be extended to
more general groups.  Other related results, on spaces of differential
forms, in the spirit of this
discussion are \cite{MPR1, MPR2} and \cite{PR,PR2}.
\ms 

Thus, a common theme of this circle of ideas is that the underlying
manifold, Riemannian or sub-Riemannian, and the collection of vector
fields $\bX$ satisfying H\"ormander's condition and defining the
corresponding sub-Laplacian, determine a metric structure and 
the most efficient way to describe smoothness of functions and
regularity of canonical operators is via scale of spaces that are
modeled by the sub-Laplacian, hence by $\bX$, itself.
\ms

In the setting of  Carnot--Carath\'eodory 
groups, and more in general of Lie groups of polynomial growth, 
endowed with the sub-Riemannian structure induced by a family $\bX$
of vector fields satisfying H\"ormander's condition, then a
Mihlin--H\"ormander multiplier theorem holds. 
This fact allowed G.\ Furioli, C.\ Melzi and A.\ Veneruso~\cite{FMV}
to introduce Besov spaces on such groups, which were later studied by
I.\ Gallagher and Y.\ Sire~\cite{GS}.  The theory was recently
extended to  any
unimodular Lie group by J.\ Feneuil~\cite{Feneuil}.  

This work aims to contribute to the analysis of function spaces on
general noncompact Lie groups, hence including the nonunimodular groups, with
Haar measures of exponential growth. 

Concerning the function spaces, their algebra properties are of great
importance, in particular in application to well-posedness results for 
nonlinear differential equations.  In this direction, a remarkable 
paper is~\cite{CRTN} by T.\ Coulhon, E.\ Russ and V.\ Tardivel-Nachev,
where they proved algebra properties for the Sobolev spaces, in
particular
on any
unimodular Lie group.    The algebra properties were
extended to the scale of Besov spaces on groups of polynomial growth
in~\cite{GS} and in~\cite{Feneuil} on
unimodular Lie groups.

A number of the aforementioned  results were also obtained in the context of doubling measure metric spaces with the reverse doubling property, see e.g.~\cite{HMY, MY} and in the
setting of Riemannian 
manifolds of bounded geometry, see
e.g.~\cite{Triebel1,Triebel2,Triebel3},~\cite{CRTN}, and references
therein.  
On the other hand, not much is
known in the setting of a sub-Riemannian manifold.  
This work is part of  a program 
\cite{PV}, \cite{BPTV} and \cite{BPV}, whose  main long term goal is
to address this type
of questions on a sub-Riemannian manifold, and we started 
with the case of a
general Lie group. The paper~\cite{PV} studies Sobolev spaces with respect to the
sum of squares sub-Laplacian,  results then extended to Sobolev spaces
with respect to sub-Laplacians with drift in~\cite{BPTV}, while
in~\cite{BPV} we develop the theory of Besov and Triebel--Lizorkin
spaces with respect to sub-Laplacians with drift, that we further
analyse in this work.

We conclude this part of the introduction by pointing out that the literature in
this area is extremely vast, and it is just impossible to give credit
to all the authors that have contributed to its development.  We
apologise to everyone whom we did not explicitly mention.
\ms

Let $G$ be a noncompact connected Lie group  and let
$\bX= \{X_1, \dots ,X_\ell\}$ be a family of linearly
independent left-invariant vector fields on $G$ satisfying
H\"ormander's condition. We denote by $\delta$ the modular function on
$G$.  Let $\rho$ be a right Haar measure of
$G$, let $\chi$ be a continuous positive character of $G$, and
consider the measure $\mu_\chi$ defined by the relation
$d \mu_\chi = \chi d \rho$. 
 Consider now
the differential operator 
\begin{equation}\label{Deltachi}
\Dc = - \sum_{j=1}^\ell (X_j^2 + c_jX_j) \,,
\end{equation}
with domain 
$C_c^\infty(G)$, where 
$c_j=(X_j \chi)(e)$, $j=1,\dots, \ell$, and  $e$ is the identity
of $G$.

 This operator was introduced by W.\ Hebisch, G.\ Mauceri and
S.\ Meda in~\cite{HMM}, where they showed that $\Dc$ is essentially
self-adjoint on $L^2(\mu_\chi)$.  Moreover, they proved that if a
sub-Laplacian 
with drift is symmetric on $L^2(\mu)$ for a positive measure $\mu$ on
$G$, 
then necessarily $\mu=\mu_\chi$ for a positive character
$\chi$ on $G$,  and moreover the drift has the form 
$X:=\sum_{j=1}^\ell c_jX_j$, where 
$c_j=(X_j \chi)(e)$, $j=1,\dots, \ell$, as in \eqref{Deltachi}. 
Notice that when the character
$\chi$ is the modular function, $\mu_\delta= \lambda$ is a left Haar
measure and the operator $\Delta_\delta$ coincides with 
the intrinsic hypoelliptic Laplacian
associated with the Carnot-Carath\'eodory metric induced on $G$ by the
vector fields $\bX$ --- see \cite{ABGR}.  The operator $\Delta_\delta$ 
is the natural
substitute of the Laplacian on a general Lie group $G$. 
This also reflects on the
fact that the measure $\lambda$ is privileged among the measures
$\mu_\chi$. As shown in~\cite{ABGR, BPTV}, $\Delta_\delta$ is not a
sum-of-squares operator unless the group is unimodular.   In this
paper we continue the study of function spaces associated with $\Dc$
for a generic positive continuous character $\chi$.  
The more general treatment allows one extra flexibility, 
see e.g.\ the embedding results, Theorems 1.1 and 4.4 
in~\cite{BPTV} and Theorems 5.2 and 5.3 in~\cite{BPV}, and at the same
time, 
highlights the naturality of $\Delta_\delta$.\ms

In this paper we further develop the investigation of Besov and Triebel--Lizorkin
spaces on $G$, defined in terms of $\Dc$, spaces that were
introduced  by the authors in the recent paper \cite{BPV}.

We
prove characterizations of the norms in term of finite differences
(Theorems \ref{Slocqalpha-charac} and
\ref{finite-difference-charac-Besov}), the density of the test functions in Besov and Triebel--Lizorkin spaces, 
and the boundedness
of  a simplified version of the 
local Riesz transforms (Theorem \ref{iso-prop}).

The plan of the paper is as follows.
In the next section we recall the basic facts about our setting and in
particular the heat semigroup
generated by $\Dc$. In Section \ref{BTL-spaces} we recall the
definitions of Besov and 
Triebel--Lizorkin spaces, and the results of \cite{BPV} needed in the
present work. 
In Section \ref{fin-diff-sec} we prove finite difference
characterizations for the Besov and Triebel--Lizorkin spaces. Such
characterizations are then used in 
Section \ref{density-sec} to show that test functions are dense in such
spaces, and in Section \ref{iso-sec} we prove an isomorphism 
result  and the
boundedness of the aforementioned version of the local Riesz transforms
for both scales of Besov and Triebel--Lizorkin spaces.
We conclude by mentioning some directions for future work.
\ms

{\em Foreword by the second named author.}
Soon after getting my  Ph.~D., I obtained a position at the Politecnico
in Torino, where Fulvio had been for a number of years.  He was my
main reason for seeking this position at the Politecnico.  I immediately
found myself immersed in a very pleasant enviroment, with Fulvio being
the organiser of many activities, such as  advanced
courses, regular seminars, and the visits of many leading mathematicians.  
I was exposed to a
flurry of recent and as well as ongoing research, on a variety of different
topics.  This gave me the possibility of meeting and interacting with many 
experts.   Fulvio personally introduced me to this
world, taking the time to explain to me a lot of mathematics, while
advising and guiding me. I have always
been very impressed by his
poise, kindness, and, most of all, generosity in teaching all the younger
mathematicians who had the fortune to
interact with him. 
He has had a great
impact
 on me, both 
 professionally and personally.

I wish to express to Fulvio my most sincere gratitude for all
he has taught me, and 
for his invaluable friendship.

\ms

\section{Basic facts and definitions}\label{BFDs-sect}
\ms 

Let $G$ be any Lie group with identity element $e$.  
We denote by $\rho$  a
right  Haar measure, and by $\delta$ the modular function.  We let
$\lambda$ be the left Haar measure such that
$d\lambda= \delta d\rho$.  We recall that $\delta$
is a smooth positive character, that is, a smooth group homomorphism of
$G$ onto $\bbR^+$. If $\chi$ is any continuous positive character of
$G$, then $\chi$ is automatically smooth.  For any such $\chi$, we
define $\mu_\chi$ to be the measure whose density with respect to
$\rho$ is $\chi$, that is, $d\mu_\chi=\chi d\rho$.  Notice that
$\mu_1=\rho$ and $\mu_\delta=\lambda$.  

We fix once for all a family of left-invariant linearly independent
vector fields 
$\bX= \{X_1,\dots,X_\ell\}$
satisfying H\"ormander's condition. 
These vector fields induce the Carnot--Charath\'eodory distance, 
denoted by $d_C$, 
which turns out to be
 left-invariant.
Then, for $x\in G$ we set $|x|=d_C(x,e)$ and we denote with $B(x,r)$ the 
ball centered at $x$ and of radius $r>0$.  If $x=e$ and $r>0$, we write
$B_r=B(e,r)$, and define $V(r)=\rho(B_r)$.  In general, we denote by
$B(c_B,r_B)$ the ball with center $c_B$ and radius $r_B$, in the
metric $d_C$.

It is known that there exist two
constants $d,D>0$ such that
\begin{equation}\label{d-D}
\begin{cases}
\rho(B_r)\approx r^d \quad & \text{if } r\in (0,1] \cr
\rho(B_r)\lesssim e^{Dr} &  \text{if } r\in (1,+\infty) \,,
\end{cases} 
\end{equation}
see \cite{Gui, Var}. It is worth to point out that $d=d(\bX,G)$, while $D=D(G)$.

For any quantities $A$ and $B$, we  write $A\lesssim B$ to indicate
that there exists a constant $c>0$, independent of the relevant parameters,
such that $A\leq c \,B$. If
$A\lesssim B$ and $B\lesssim A$, we write $A\approx B$.  In order to
emphasize the dependence on a given parameter, say $R$, we write
$\lesssim_R$, and analogously for the other cases.

We observe that, having fixed $R>0$,  every character $\chi$ satisfies
the estimates  
\begin{equation}\label{equiv-char}
\chi(x) \approx_R\chi (y) 
\end{equation}
for all
$x,y\in G$ such that 
$d_C(x,y)\leq R$.  This equivalence easily implies that $(G,d_C,\mu_\chi)$ is locally doubling, that is, for all $0<r<R$ and $x_0\in G$,
\begin{equation}\label{local-doubling-mu-chi}
\mu_\chi (B(x_0,2r)\lesssim_R \mu_\chi(B(x_0,r))\,.
\end{equation}

Having fixed $\bX$, we consider the operator $\Dc$ defined
in~\eqref{Deltachi}.  With an abuse of notation, we still denote 
by $\Dc$ its smallest closed  extension on
$L^p(\mu_\chi)$, where, for  $p\in(1,+\infty)$, $L^p(\mu_\chi)$ denotes
  the standard Lebesgue space.  The space $L^\infty$ is the space of
  $\rho$-essentially bounded functions.  We refer to~\cite{HMM} and
  \cite{BPTV} for further details about $\Dc$.

We set $\sJ=\{1,\dots,\ell\}$ and we say that a multi-index
$J=(j_1,\dots,j_m)\in \sJ^m$ if $j_k\in \sJ$ for
$k=1,\dots,m$. Moreover, we write
$$
X_J=X_{j_1}\cdots,X_{j_m} \,.
$$

Next, we observe that, since $\Dc$ is left-invariant, the associated heat
semigroup admits a convolution kernel  $p_t^\chi\in \mathcal{D'}(G)$,
i.e.\   
$$
e^{-t\Dc} f = f* p_t^\chi \,.
$$
It is known (see~\cite{BPTV}) that
\begin{equation}\label{ptgamma_pt}
p_t^\chi= e^{-c_X t/4}  \chi^{-1/2} p_t \,,
\end{equation}
where $c_X = \big(\sum_{j=1}^{\ell}(X_j\chi(e))^2\big)^{1/2}$, and
$p_t$ is the convolution kernel of the heat semigroup in the case
$\chi=1$, for which the estimates in~\cite{VSCC} are available. 

We recall the expressions of the convolution on $G$.  We have 
\begin{align}
f*g(x)
& = \int_G f(xy^{-1}) g(y)\, d\rho(y) 
= \int_G f(y^{-1}) g(yx)\, d\rho(y) \notag \\
& = \int_G f(y) g(y^{-1}x)\, d\l(y) = \int_G f(xy) g(y^{-1})\, d\l(y)
\,.  \ms
\label{convolutions}
\end{align}

The following result is essentially Lemma 3.1 in~\cite{BPV}. 
\begin{lem}\label{estimates_pt-chi} 
The following properties hold:
\begin{itemize}
\item[(i)] $(e^{-t\Dc})_{t>0}$ is a diffusion semigroup on $(G, \mu_\chi)$;
\smallskip
\item[(ii)] for every $r>0$, $\sup_{B_r}\chi=e^{c_X r}$;
\item[(iii)] there exist two constants $c_1, c_2>0$ such that 
$$
(\delta \chi^{-1})^{1/2}(x)\, V(\sqrt t)^{-1}e^{-c_1|x|^2/t} 
\lesssim p_t^\chi (x)\lesssim (\delta \chi^{-1})^{1/2}(x) \, V(\sqrt t)^{-1}e^{-c_2|x|^2/t}
$$
 for every $t\in (0,1)$ and $x\in G$;
\item[(iv)]  given $m\in \bbN$, there exist a positive constant $b=b_m$ such that for $ x\in G,\, J\in \sJ^m$
$$
|X_J p_t^\chi (x)| \lesssim (\delta\chi^{-1})^{1/2}(x) t^{-m/2}\,V(\sqrt{t})^{-1}
\,e^{-b|x|^2/t}  \qquad \text{for } t\in(0,1) \,;
$$ 
\item[(v)] there exists  $c_3>0$ such that  for every $t\in (0,1)$ and $x\in G$,
$$
\Big|\frac{\p }{\p t}\, p_t^\chi(x)\Big| \lesssim
(\delta \chi^{-1})^{1/2}(x) t^{-1} V(\sqrt{t})^{-1}  e^{-c_3 |x|^2/t}
\,.
$$
\end{itemize}
\end{lem}

Only the bound in (v) needs to be justified.  It follows from the
estimates 
 for the unweighted heat kernel provided in~\cite[Section IX]{VSCC}
 (see also~\cite[Section 2.1, estimates (ii)-(iv)]{PV}) together
 with~\cite[(2.8)]{BPTV}.
\ms

\begin{defn}{\rm
We define the Schwartz space $\cS(G)$ as the space of functions
$\vp\in C^\infty(G)$ such that for all $n,m\in\bbN$, $J\in\sI^m$  the seminorms
$$
\cN_{J,n}(\vp) = \sup_{x\in G} e^{n|x|} |X_J \vp(x)| 
$$
are finite. The space $\cS'(G)$ is defined as the dual space of
$\cS(G)$.
}
\end{defn}

As a consequence of the gaussian estimates (iv) in Lemma
\ref{estimates_pt-chi}, we have the following simple lemma.

\begin{lem}\label{pt-Schwartz-cont}
For all $t>0$,
$p_t^\chi\in\cS$.  Moreover,
$e^{-t\Dc}:\cS\to\cS$ is bounded, with 
seminorms uniformly bounded for $t\in[\eps,R]$, for any $0<\eps<R$.   Therefore, 
$e^{-t\Dc}$ extends to a continous map
$e^{-t\Dc}:\cS'\to\cS'$,
for all $t>0$.
\end{lem}

\proof
We indicate the argument for sake of completeness.  Given $n\in\bbN$
and a multi-index $J$,  we have 
\begin{align}
e^{n|x|} | X_J (\vp*p^\chi_t)(x)|
& = e^{n|x|} |( \vp*X_J p^\chi_t)(x)| \notag \\
& \le  \int_G
e^{n|xy^{-1}|} |\vp(xy^{-1})| e^{n|y|} |X_J p^\chi_t(y)| \, d\rho(y) \notag\\
& \le \cN_{0,n}(\vp)  \int_G e^{n|y|} |X_Jp^\chi_t(y)|\, d\rho(y) \notag \\
& \lesssim t^{-|J|/2} \cN_{0,n}(\vp)\,, \label{Sch-est}
\end{align} 
where the last inequality is obtained arguing as in the proof of Lemma
3.4 in \cite{BPV}.  The conclusions now follow easily.  
\epf

\begin{defn}\label{W-t-m}{\rm
For $m\in\bbN$ and $t>0$, we define the operator $W_t^{(m)}$ by
setting
$$
W_t^{(m)}
=
(t\Dc)^m e^{-t\Dc}\,.
$$   
}
\end{defn}
For $m\in\bbN$ and $t>0$, 
$W_t^{(m)}:\cS\to\cS$ is bounded, and therefore it extends to a
continuous map $W_t^{(m)}:\cS'\to\cS'$.  We also observe that, for
$f\in\cS'$, $W_t^{(m)}f$ is a $C^\infty$ function, for all $t>0$ and
$m\in\bbN$. \ms

We also recall the definition of the Littlewood--Paley--Stein
$g$-function.
Given a positive integer $k$, for $f\in\cS$ we set
\begin{equation}\label{def-g}
g_k(f) =
\bigg( \int_0^{+\infty} 
\big| W^{(k)}_s  f\big|^2\, \frac{ds}{s} \bigg)^{1/2}
 \,.
\end{equation}
Since $\Dc$ generates a  symmetric diffusion
semigroup, for $p\in(1,+\infty)$, 
$g_k$  satisfies the estimate
\begin{equation}\label{g-fct-est}
\| g_k(f)\|_{L^p(\mu_\chi)} \approx \| f\|_{L^p(\mu_\chi)} \,,
\end{equation}
see~\cite{Stein}, and also \cite{Meda}. 

\ms

\section{Triebel--Lizorkin and Besov spaces on $G$}\label{BTL-spaces}
\ms

Here and in what follows, given a measure space $(\Omega,\nu)$ and a
Banach space $\cX$, for $p\in[1,+\infty]$, 
we denote by $L^p\big( (\Omega,\cX),\nu\big)$ the
space of measurable functions $f:\Omega\to\cX$ such that
$$
\| f\|_{L^p ( (\Omega,\cX),\nu)} :=
\Big\{ \int_\Omega \| f(\omega)\|_{\cX}^p \, d\nu(\omega) \Big\}^{1/p}
<\infty \,
$$
when $p\in[1,+\infty)$,
with the obvious modification if $p=+\infty$.  We also 
denote by $[\tau]$ the integral part of $\tau\ge0$.
\ms

We are now in the position to introduce the Triebel--Lizorkin and Besov
spaces on $G$, defined in terms of the sub-Laplacian $\Dc$,
see \cite{BPV}.  

\begin{defn}\label{def-TL-B}{\rm
Let $p,q\in[1,+\infty]$, and $\alpha\ge0$.  Then we define:\smallskip

\noindent
{\rm (1)} the Triebel--Lizorkin space $F^{p,q}_\alpha(\mu_\chi)$ as
$$
F^{p,q}_\alpha (\mu_\chi)= \Big\{ f\in\cS'(G):\,  t^{-\alpha/2}
W_t^{([\alpha/2]+1)} f \in L^p \Big( \big(G, L^q
\big((0,1),dt/t\big) ,\,
\mu_\chi \Big)\ \text{and}\ e^{-\frac12\Dc}f \in L^p(\mu_\chi) \Big\}
$$
endowed with the norm
\begin{equation}\label{TL-norm}
\| f\|_{F^{p,q}_\alpha} := 
\sF^{p,q}_{\alpha} (f) +
\| e^{-\frac12\Dc} f\|_{L^p(\mu_\chi) } \,,
\end{equation}
where
$$
\sF^{p,q}_{\alpha} (f) := 
\Big\| \bigg( \int_0^1 \Big( t^{-\alpha/2}  \big|
W_t^{([\alpha/2]+1)}  f
\big| \bigg)^q \, \frac{d t}{t}\bigg)^{1/q} \Big\|_{L^p(\mu_\chi)}
$$
if $q<+\infty$, while 
$$
\sF^{p,\infty}_{\alpha} (f) := \big\| \sup_{t\in (0,1)} t^{-\alpha/2} |
W_t^{([\alpha/2]+1)}  f|\big\|_{L^p(\mu_\chi)} \,;
$$

\noindent
{\rm (2)} the Besov space $B^{p,q}_\alpha(\mu_\chi)$ as
$$
B^{p,q}_\alpha (\mu_\chi)= \Big\{ f\in\cS'(G):\,  t^{-\alpha/2}
W_t^{([\alpha/2]+1)} f \in L^q \Big(  \big( (0,1), L^p(G,\mu_\chi)
\big), \, dt/t\Big)\ \text{and}\ e^{-\frac12 \Dc}f \in L^p(\mu_\chi) \Big\}
$$
endowed with the norm
\begin{equation}\label{Besov-norm}
\| f\|_{B^{p,q}_\alpha} := 
\sB^{p,q}_{\alpha} (f)  +
\| e^{-\frac12 \Dc} f\|_{L^p(\mu_\chi) } \,,
\end{equation}
where
$$
\sB^{p,q}_{\alpha} (f) :=
\left( \int_0^1 \left(t^{-\alpha/2}\,  \| W_t^{([\alpha/2]+1)}
    f\|_{L^p(\mu_\chi)}\right)^q \, \frac{dt}{t}\right)^{1/q} 
$$
if $q<+\infty$, while 
$$
\sB^{p,\infty}_{\alpha} (f) :=\sup_{t\in (0,1)}
t^{-\alpha/2} \,\|  W_t^{([\alpha/2]+1)}   f\|_{L^p(\mu_\chi)} \,.
$$
}
\end{defn}

We emphasize that, when $p\in(1,+\infty)$ and $\alpha\ge0$,
the Triebel--Lizorkin space $F^{p,2}_\alpha(\mu_\chi)$ coincides with
the Sobolev space $L^p_\alpha(\mu_\chi)$ defined in~\cite{BPTV}, with
equivalence of norms, see~\cite[Theorem 5.2]{BPV}.
\ms

We now recall the main results in~\cite{BPV} about equivalence of
norms in Besov and Triebel--Lizorkin spaces.  The following is~\cite[Theorem
4.1]{BPV}.
\begin{thm} \label{teo-equiv1}
Let $\alpha> 0$, $m>\alpha/2$ be an integer, $t_0\in [0,1)$ and $q\in [1,+\infty]$.
\begin{itemize}
\item[(i)] If $p\in (1,+\infty)$, then the norm $\|f\|_{F_\alpha^{p,q}} $ is equivalent to the norm
\begin{equation}\label{TL1}
\Big\| \bigg( \int_0^1 \Big( t^{-\alpha/2}| W_t^{(m)}
    f|\Big)^q \, \frac{d t}{t}\bigg)^{1/q}\Big\|_{L^p(\mu_\chi)}
+ \|e ^{-t_0\Dc}f\|_{L^p(\mu_\chi)}.
\end{equation}
\item[(ii)] If $p \in [1,+\infty]$, then the norm $\|f\|_{B_\alpha^{p,q}} $ is equivalent to the norm
\begin{equation}\label{Besov2}
\bigg( \int_0^1 \Big( t^{-\alpha/2} \| W_t^{(m)}
f\|_{L^p(\mu_\chi)}\Big)^q \, 
\frac{d t}{t}\bigg)^{1/q} + \|e ^{-t_0\Dc}f\|_{L^p(\mu_\chi)}.
\end{equation}
\end{itemize}
If $\alpha=0$, the norms $\|f\|_{F_\alpha^{p,q}} $ and
$\|f\|_{B_\alpha^{p,q}} $ are equivalent to those
in~\eqref{TL1} and \eqref{Besov2} respectively  provided  $t_0\in (0,1)$.
\end{thm}

The next result concerns a discretization of the norm that resembles 
the Littlewood--Paley characterization of Besov and Triebel--Lizorkin
spaces in the classical cases.
In our case, for $j\in\bbN$, the operators $W^{(m)}_{2^{-j}}$ play the
role of the operators $\bigtriangleup_j$ in the classical
Littlewood--Paley decomposition, while $e^{-t_0\Dc}$ plays the
role of $S_0$; see e.g.~\cite{Grafakos} for such notation in the case
of $\bbR^d$.

We point out that in the case of $\Dc$ the
classical Littlewood--Paley characterization of Besov and Triebel--Lizorkin
spaces cannot hold since any bounded spectral multiplier of $\Dc$ on
$L^p(\mu_\chi)$, with $p\neq2$, admits a
holomorphic extension to a parabolic region in $\bbC$, see~\cite{HMM}.  This is~\cite[Theorem
4.2]{BPV}.

\begin{thm} \label{teo-equiv2}
Let $\alpha>  0$, $m>\alpha/2$ be an integer, $t_0\in [0,1)$ and $q\in [1,+\infty]$.
\begin{itemize}
\item[(i)]  If $p\in (1, +\infty)$, then the norm $\|f\|_{F_\alpha^{p,q}} $ is equivalent to the norm
\begin{equation}\label{TL2}
\Big\| \bigg(\sum_{j=0}^\infty \Big( 2^{j \alpha/2}
|W^{(m)}_{2^{-j}} f|\Big)^q \bigg)^{1/q} \Big\|_{L^p(\mu_\chi)}
+ \| e^{-t_0\Dc}f\|_{L^p(\mu_\chi)} \,.
\end{equation}
\item[(ii)] If $p\in [1,+\infty]$, then the norm $\|f\|_{B_\alpha^{p,q}} $ is equivalent to the  norm
\begin{equation}\label{Besov3}
\bigg(\sum_{j=0}^\infty \bigg( 2^{j \alpha/2}
\|W^{(m)}_{2^{-j}} f\|_{L^p(\mu_\chi)} \bigg)^q \bigg)^{1/q} +
 \| e^{-t_0\Dc}f\|_{L^p(\mu_\chi)} \,.
\end{equation}
\end{itemize}
If $\alpha=0$, the norms 
$\|f\|_{F_\alpha^{p,q}} $ and
$\|f\|_{B_\alpha^{p,q}} $ 
are equivalent respectively to those in~\eqref{TL2} 
and~\eqref{Besov3} provided $t_0\in (0,1)$.
\end{thm}

\section{Finite differences characterizations} \label{fin-diff-sec} 
\ms

In this section we prove characterizations for the spaces
$F^{p,q}_\alpha$ and $B^{p.q}_\alpha$ in terms of finite differences. 
Such characterizations provide a key tool for the proof of the density
lemma of Section \ref{density-sec}.   We begin by introducing the
finite difference operator.\ms

Given a measurable function $f$, for $x,y\in G$ we define
\begin{equation}\label{diff-op-eq}
\oD_yf(x)= f(xy^{-1})-f(x) \,.
\end{equation}
\ms

\subsection{Characterization of Triebel--Lizorkin norm by differences }
For $q\in [1,\infty]$ and $\alpha\in(0,1)$, we define the functional 
\begin{equation}\label{Salphaqloc}
\sS^{{\rm loc}, q}_{\alpha}f(x) =\left(\int_0^1\left[ \frac{1}{u^{\alpha}
      V(u)}\int_{|y|<u}|\oD_y f(x)|\, d\rho(y) \right]^q\,
  \frac{d u}{u}\right)^{1/q}  \,.
\end{equation}
We point out that in the case $q=2$, such functional coincides with
the classical functional $\sS^{\rm loc}_{\alpha}$ used to characterize
the Sobolev norm, see \cite{CRTN} for the unimodular case, and
\cite{BPTV} for the nonunimodular (and weighted) case.

The first result of this section is the characterization of the
Triebel--Lizorkin norm of $F^{p,q}_\alpha$ 
in terms of the $L^p(\mu_\chi)$-integrability of the functional
$\sS^{{\rm loc}, q}_{\alpha}$.  To the best of our knowledge, such
characterization is new even in the case of $\bbR^d$.

\begin{thm}\label{Slocqalpha-charac}
For every $p,q\in (1,\infty)$, $\alpha \in (0,1)$, we have
$$
\|f\|_{F^{p,q}_\alpha} \approx    \|\sS^{{\rm loc},
  q}_{\alpha}f\|_{L^p(\mu_\chi)}+\|f\|_{L^p(\mu_\chi)} \,. 
$$
\end{thm}

\proof
Set 
$$
H_{\alpha,q} f:= \left(\int_0^1 \left( t^{-\alpha/2}
    | W^{(1)}_t f| \right)^q  \, \frac{d
    t}{t}\right)^{1/q} \,,
$$
and observe that, since $\alpha\in(0,1)$, 
$$
\|f\|_{F^{p,q}_\alpha} \lesssim \| H_{\alpha,q}f \|_{L^p(\mu_\chi)}
+ \|f\|_{L^p(\mu_\chi)} \,.
$$ 

{\em Step 1.} We shall prove that  for all $f\in
F^{p,q}_\alpha$, 
\begin{equation}\label{claim1}
\|f\|_{F^{p,q}_\alpha}
\lesssim   \|\sS^{{\rm loc},
  q}_{\alpha}f\|_{L^p(\mu_\chi)}+\|f\|_{L^p(\mu_\chi)} \,, 
\end{equation}
by showing that
\begin{equation}\label{step-1}
\|H_{\alpha,q} f\|_{L^p(\mu_\chi)}\lesssim \|\sS^{{\rm loc},
  q}_{\alpha}f\|_{L^p(\mu_\chi)}+\|f\|_{L^p(\mu_\chi)}\,.
\end{equation}

We first notice that for every $t\in (0,1)$ and $x\in G$, since
$\frac{\p}{\p t}\int_Gp_t^{\chi} \, d\rho =0$, we have  
\begin{align*} 
| \Dc e^{-t\Dc}f(x)| 
& = \Big| \frac{\p}{\p t}e^{-t\Dc}f(x) \Big|
= \Big| \frac{\p}{\p t} \left(  \int_G f(xy^{-1})
    p_t^{\chi}(y)\, d \rho(y)- \int_G f(x) p_t^{\chi}(y) \, d
    \rho(y)  \right) \Big| \\  
& \leq \int_G|\oD_yf(x)|  \Big|\frac{\p p_t^{\chi}(y)}{\p t}\Big|\, d \rho(y). 
\end{align*}
Using the estimates (v) of Lemma~\ref{estimates_pt-chi} we have 
\begin{align*} 
\big( H_{\alpha,q}f(x) \big)^q
&\lesssim \int_0^1 t^{-q\alpha/2}  V(\sqrt t)^{-q}
\bigg(\int_{|y|<\sqrt t}|\oD_yf(x)| (\delta \chi^{-1})^{1/2}(y)
  e^{-c_4|y|^2/t} d \rho(y)\bigg)^q \,\frac{d t}{t}
\\   &\quad 
+\sum _{k=0}^{\infty} \int_0^1 t^{-q\alpha/2} V(\sqrt 
t)^{-q}\bigg( \int_{ 2^k \sqrt{t} <|y|<2^{k+1}\sqrt t }
  |\oD_yf(x)|(\delta \chi^{-1})^{1/2}(y)   \,  e^{-c_4|y|^2/t}
  d \rho(y)\bigg)^q \,\frac{d t}{t}\\ 
&\lesssim  \int_0^1 t^{-q\alpha/2}  V(\sqrt
t)^{-q}\bigg(\int_{|y|<\sqrt t}|\oD_yf(x)|   d
  \rho(y)\bigg)^q \,\frac{d t}{t}\\  
&\quad +\sum _{k=0}^{\infty}e^{-c_4 2^{2k}}\int_0^1
t^{-q\alpha/2} V(\sqrt t)^{-q}\bigg(\int_{|y|<2^{k+1}\sqrt
    t}|\oD_yf(x)| (\delta \chi^{-1})^{1/2}(y)d
  \rho(y)\bigg)^q\,\frac{d t}{t}. 
\end{align*}
By the change of variables $u=2^{k+1}\sqrt t$ we obtain
\begin{align*}
\big( H_{\alpha,q}f(x) \big)^q 
&\lesssim   \int_0^1 \frac{1}{u^{q\alpha}V(u)^q}\bigg(\int_{|y|<
    u}|\oD_yf(x)|  d \rho(y)\bigg)^q\, \frac{d u}{u}
\\  &\quad 
+\sum
_{k=0}^{\infty}e^{-c_4 2^{2k}}\int_0^{2^{k+1}}\frac{2^{(k+1)q\alpha}}{
u^{q\alpha}V(2^{-k-1}u)^q  }\bigg(\int_{|y|<u}|\oD_yf(x)|
  (\delta \chi^{-1} )^{1/2}(y)d \rho(y)\bigg)^q \,  \frac{d
  u}{u}\\  
&\lesssim \big( \sS^{{\rm loc}, q}_{\alpha}f(x) \big)^q
+\sum _{k=0}^{\infty}e^{-c_4 2^{2k}}\int_0^{1}
\frac{2^{(k+1)q\alpha}}{ u^{q\alpha}V(2^{-k-1}u)^q  }
\bigg(\int_{|y|<u}|\oD_yf(x)| \, d \rho(y)\bigg)^q\,
\frac{d u}{u}\\  
&\quad +\sum _{k=0}^{\infty}e^{-c_4 2^{2k}}\int_1^{2^{k+1}}
\frac{2^{(k+1)q\alpha}}{u^{q\alpha}V(2^{-k-1}u)^q  }
\bigg(\int_{|y|<u} |\oD_yf(x)|  (\delta \chi^{-1} )^{1/2}(y)d
  \rho(y)\bigg)^q \,  \frac{d u}{u}. 
\end{align*} 
By the estimates~\eqref{d-D}, we obtain that
\begin{align}
\big( H_{\alpha,q}f(x) \big)^q 
&\lesssim \big( \sS^{{\rm loc}, q}_{\alpha}f(x)\big)^q +\big( \sS^{{\rm
    loc}, q}_{\alpha}f(x) \big)^q
\sum _{k=0}^{\infty}e^{-c_4 2^{2k}}2^{(k+1)(q\alpha+qd)} \notag \\  
&\quad +\sum _{k=0}^{\infty}e^{-c_4 2^{2k}}   2^{(k+1)(q\alpha+qd)}
\int_1^{2^{k+1}}  \frac{1}{  u^{q\alpha+qd}  }  |f(x)|^q
\bigg(\int_{|y|<u}  (\delta \chi^{-1})^{1/2}(y)\, d \rho(y)\bigg)^q
\, \frac{d u}{u} \notag \\  
&\quad +\sum _{k=0}^{\infty}e^{-c_4 2^{2k}}   
2^{(k+1)(q\alpha+qd)} \int_1^{2^{k+1}}  \frac{1}{  u^{q\alpha+qd}  }
\bigg(\int_{|y|<u}|f(xy^{-1})|  (\delta \chi^{-1})^{1/2}(y)\, d
  \rho(y)\bigg)^q \, \frac{d u}{u} \notag \\  
&\lesssim \big( \sS^{{\rm loc}, q}_{\alpha}f(x) \big)^q +\sum_{k=0}^{\infty}J_k(x)+\sum
_{k=0}^{\infty}I_k(x) \,. \label{dec-G-alpha-q}
\end{align}
By the growth estimates of characters in Lemma~\ref{estimates_pt-chi}
(ii), we deduce that there exists $C>0$ such that
\begin{align} 
\sum _{k=0}^{\infty} J_k(x)
&\lesssim  \sum _{k=0}^{\infty} e^{-c_4 2^{2k}}   2^{(k+1)(q\alpha+qd)}
\int_1^{2^{k+1}}  \frac{1}{  u^{q\alpha+qd}  }  |f(x)|^q   u^{qd}e^{qC
  u} \frac{d u}{u} \notag \\ 
&\lesssim \sum _{k=0}^{\infty} |f(x)|^q   e^{-c_4 2^{2k}}
2^{(k+1)(q\alpha+qd)}  e^{qC 2^k} \notag \\
& \lesssim |f(x)|^q \,.\label{Jk}
\end{align}
We now notice that there exists $c >0$ such that 
$$
\Big\| \bigg( \sum_{k=0}^{\infty}I_k\bigg)^{1/q}\Big\|_{L^p(\mu_\chi)}
\lesssim \sum_{k=0}^{\infty} e^{-c 2^{2k}}\Big\| \bigg(
    \int_1^{2^{k+1}}   \bigg(\int_{|y|<u}|f(\cdot\, y^{-1})|
      (\delta \chi^{-1})^{1/2}(y)\, d \rho(y)\bigg)^q \bigg)^{1/q}
  \,  d u \Big\|_{L^p(\mu_\chi)}. 
$$
Here and in the rest of this work, we denote by  ${\mathbf 1}_E$  
the characteristic
function of the measurable set  $E$.
For every integer $k$, by Minkowski's integral inequality, we get 
\begin{align*}
&\bigg(   \int_1^{2^{k+1}} \bigg(\int_G |f(x y^{-1})| {\mathbf 1}_{B_u}(y)
    (\delta \chi^{-1})^{1/2}(y)\, d \rho(y)\bigg)^q  \, d u
\bigg)^{1/q} \\ 
&\qquad \lesssim \int_G \bigg(
  \int_1^{2^{k+1}}|f(xy^{-1})|^q {\mathbf 1}_{B_u}(y)(\delta \chi^{-1})^{q/2}(y)\,
  d u \bigg)^{1/q}\, d\rho(y)\\ 
&\qquad\lesssim \int_{B_1} |f(xy^{-1})| (\delta \chi^{-1})^{1/2}(y)
\bigg(  \int_1^{2^{k+1}} \, du \bigg)^{1/q} \, d\rho(y)\\ 
&\qquad\qquad +\int_{1<  |y|< 
  2^{k+1}} |f(xy^{-1})| (\delta \chi^{-1})^{1/2}(y)  \bigg(
  \int_{|y|}^{2^{k+1}}\, du \bigg)^{1/q}\, d\rho(y)\\ 
&\qquad\lesssim 2^{k/q} \int_{B_1}|f(xz)| \, d\l(z) + 
2^{k/q} \int_{1<  |z|<  2^{k+1}} |f(xz)| (\delta^{-1}\chi)^{1/2}(z)\,
d\l(z)\,. 
\end{align*}
By applying again Minkowski's inequality, we then obtain that
\begin{align*}
 &\Big\|\bigg( \int_1^{2^{k+1}} \bigg(\int_G |f(\cdot\, y^{-1})|
       {\mathbf 1}_{B_u}(y) (\delta \chi^{-1} )^{1/2}(y)\, d\rho(y)
       \bigg)^q\,du  \bigg)^{1/q} \Big\|_{L^p(\mu_\chi)}\\ 
&\lesssim 2^{k/q}  \int_{B_1} \bigg(\int_G|f(xz)|^p
  d\mu_\chi (x)\bigg)^{1/p} \, d\l(z)+ 2^{k/q}  \int_{1<|z|<2^{k+1}}
  \bigg( \int_G|f(xz)|^p  \, d\mu_\chi (x)\bigg)^{1/p}
(\delta^{-1}\chi)^{1/2}(z)\, d\l(z)\\ 
&\lesssim
2^{k/q}\|f\|_{L^p(\mu_\chi)} + 2^{k/q}\|f\|_{L^p(\mu_\chi)}
\int_{B_{2^{k+1}}} \chi^{1/p}(\delta \chi^{-1})^{1/2}\,d\rho\\ 
&\lesssim 2^{k/q+kd} e^{C2^k} \|f\|_{L^p(\mu_\chi)} \,.
\end{align*}
We then have
\begin{equation}\label{Ik}
\begin{aligned}
\Big\| \bigg(
    \sum_{k=0}^{\infty}I_k\bigg)^{1/q}\Big\|_{L^p(\mu_\chi)}\lesssim
\sum_{k=0}^{\infty} e^{-c2^{2k}}2^{k/q+kd} e^{C2^k} \|f\|_{L^p(\mu_\chi)}
\lesssim \|f\|_{L^p(\mu_\chi)}. 
\end{aligned}
\end{equation}
In conclusion, by \eqref{Jk} and \eqref{Ik} we get~\eqref{step-1}, 
as required. 
\ms 
 
 It remains to show that for all $f\in F^{p,q}_\alpha$
\begin{equation}\label{claim2}
 \|\sS^{{\rm loc}, q}_{\alpha}f\|_{L^p(\mu_\chi)}\lesssim
 \| f\|_{F^{p,q}_\alpha}  \,.  
\end{equation}
In order to prove this, we write $f=(f-e^{-\Dc}f)+e^{-\Dc}f$ and
we estimate 
$\|\sS^{{\rm loc},  q}_{\alpha}(f-e^{-\Dc}f)\|_{L^p(\mu_\chi)}$ and 
$\|\sS^{{\rm loc}, q}_{\alpha}e^{-\Dc}f\|_{L^p(\mu_\chi)}$
separately.   \ms

{\em Step 2.} We 
 prove that 
\begin{equation}\label{step2-finite-diff-thm}
 \|\sS^{{\rm loc},  q}_{\alpha}(f-e^{-\Dc}f)\|_{L^p(\mu_\chi)}\lesssim
 \| f\|_{F^{p,q}_\alpha}  \,.  \ms
\end{equation}

Arguing similarly to~\cite[2.1.2.]{CRTN} we write 
\begin{equation}\label{m-decomposition}
f-e^{-\Dc}f 
=\sum_{m=1}^{+\infty} \int_{2^{-m}}^{2^{-m+1}}
\frac{\p}{\p t}e^{-t\Dc}f \, dt
= :\sum_{m=1}^{+\infty}
f_m\,.
\end{equation}
We then obtain
\begin{align}
\big( \sS^{{\rm loc}, q}_{\alpha}(f-e^{-\Dc}f)(x) \big)^q
&=\int_0^1  \left( \frac{1}{u^\alpha V(u)} \int_{|y|<u}|
  \oD_y(f-e^{-\Dc}f))(x)|\,
  d\rho(y)\right)^q\, \frac{d u}{u} \notag \\  
&=\sum_{j=1}^{+\infty} \int_{2^{-j}}^{2^{-j+1}}
\left( \frac{1}{u^\alpha V(u)} 
\int_{|y|<u}|
 \oD_y (f-e^{-\Dc}f)(x)|\,
  d\rho(y)\right)^q\, \frac{d u}{u} \notag \\  
&\lesssim \sum_{j=1}^{+\infty}  2^{jq\alpha}
\left( 2^{jd}  \int_{|y|<2^{-j+1}}|
  \oD_y(f-e^{-\Dc}f)(x)|\, d\rho(y)\right)^q  \notag \\ 
& 
\lesssim   \sum_{j=1}^{+\infty}   2^{jq\alpha}
\left(
\sum_{m=1}^{+\infty}   2^{jd}
\int_{|y|<2^{-j+1}}| \oD_yf_m(x)|
\,d\rho(y) \right)^q \notag \\
& 
=   \sum_{j=1}^{+\infty}   2^{jq\alpha}
\bigg( 2^{jd} \bigg( \sum_{m=1}^{2j} +
\sum_{m=2j+1}^{+\infty}   \bigg)
\int_{|y|<2^{-j+1}}| \oD_yf_m(x)|
\,d\rho(y) \bigg)^q \,,  \label{all-m}
\end{align}
where $f_m$ is defined in~\eqref{m-decomposition}.
If $m> 2j$, then 
\begin{equation}\label{m<}
 2^{jd} \int_{|y|<2^{-j+1}}| \oD_y f_m(x)| \, d\rho(y) \lesssim Mg_{m+1} (x)\,,
\end{equation}
where 
$$
 g_{m+1}=\int_{2^{-m}}^{2^{-m+1}} 
\Big|\frac{\p}{\p t}e^{-t\Dc}f \Big|\,d t \,,
$$
and 
$M$ is the local maximal function
 with respect
to the right Haar measure, 
\begin{equation}\label{MR}
M f(x)= \sup_{x\in B,\, r_B\le1} \frac{1}{\rho(B)} \int_B|f|\,
d\rho \,,
\end{equation}
which is bounded on $L^p(\mu_\chi)$ for every $p \in (1,\infty)$,
see~\cite[Subsection 5.1]{BPTV}.

In order to treat the case when $m\le 2j$, we notice that for every
$j\ge1$, $y\in B_{2^{-j-1}}$ and $x\in G$     
\begin{equation}\label{teovalormedio}
|\oD_y f_m(x)|\leq 2^{-j+1}\,\sup \big\{|X_if_m(w)|:
i=1,\dots,\ell,
\,|w^{-1}x|\leq 2^{-j+1} \big\}\,.
\end{equation}
Since, 
 \begin{equation}\label{m-decomposition-2}
f_m=\int_{2^{-m-1}}^{2^{-m}}   \frac{\p}{\p
  t}(e^{-2t\Dc}f)\, d t=2
\int_{2^{-m-1}}^{2^{-m}} e^{-t\Dc}
\frac{\p}{\p  t} (e^{-t\Dc}f)\, dt, 
\end{equation}
by applying the estimates of the heat
kernel in Lemma~\ref{estimates_pt-chi}  (v) and (iii), for every $w$ such that
$|w^{-1}x|\leq 2^{-j+1}$ we have  
\begin{align*}
|X_if_m(w)|
&\lesssim \int_{2^{-m-1}}^{2^{-m}}  \int_G \Big|
  \frac{\p}{\p t}(e^{-t\Dc}f)(z) \Big|
\big|X_ip_t^{\chi}(z^{-1}w)\big|\, d\l(z)\\ 
&\lesssim \int_G \int_{2^{-m-1}}^{2^{-m}}  \Big| 
\frac{\p}{\p    t}(e^{-t\Dc}f)(z) \Big| t^{-1/2} 
V(\sqrt t)^{-1}(\chi^{-1}\delta)^{1/2}(z^{-1}w) e^{ -b|z^{-1}w|^2/t }\, dt\,  d \l(z)\\
&\lesssim 2^{m/2}2^{md/2} \int_G
g_m(z) (\delta \chi^{-1})^{1/2}(z^{-1}x) e^{-b 2^m|z^{-1}x|^2} d\l(z)\\ 
&\lesssim 2^{m/2}e^{-c2^{-m}\Dc} g_m(x),
\end{align*}
for a suitable constant $c$. From \eqref{teovalormedio}  it follows that 
\begin{equation}\label{altrocaso} 
2^{jd} \int_{|y|<2^{-j+1}}|\oD_y f_m(x)| \, d\rho(y)
 \lesssim 2^{-j+m/2}e^{-c2^{-m}\Dc}g_m(x)\,.
\end{equation}

Thus, putting together \eqref{all-m}, \eqref{m<} and \eqref{altrocaso} we obtain
\begin{align}
& \big( \sS^{{\rm loc}, q}_{\alpha}(f-e^{-\Dc}f)(x) \big)^q
\notag \\
& \quad\lesssim
\sum_{j=1}^{+\infty} 2^{jq\alpha} \bigg(
\sum_{m=1}^{2j} 2^{-j+m/2}e^{-c2^{-m}\Dc}g_m(x)
+ \sum_{m=2j+1}^{+\infty} Mg_{m+1}(x) 
  \bigg)^q  \notag \\
& \quad
\lesssim
\sum_{j=1}^{+\infty}  2^{jq\alpha} \bigg( 
\sum_{m=1}^{2j} 2^{-j+m/2}e^{-c2^{-m}\Dc}g_m(x) \bigg)^q 
+  \sum_{j=1}^{+\infty}  2^{jq\alpha} \bigg( \sum_{m=2j+1}^{+\infty} Mg_{m+1}(x) 
  \bigg)^q \notag \\
& \quad =: \Sigma_1(x) + \Sigma_2(x) \,. \label{CRTN-argument}
\end{align}
Now we apply H\"older's inequality to see that, for any $\eps>0$
\begin{align*}
\bigg( 
\sum_{m=1}^{2j} 2^{-j+m/2}e^{-c2^{-m}\Dc}g_m(x) \bigg)^q 
& \le  \bigg(  \sum_{m=1}^{2j} 2^{\eps mq'} \bigg)^{q/q'} 
 \sum_{m=1}^{2j} 2^{-jq+mq/2 -\eps mq}\big(
 e^{-c2^{-m}\Dc}g_m\big (x) \big)^q \\
& \lesssim \sum_{m=1}^{2j} 2^{2j\eps q -jq+mq/2 -\eps mq}\big(
 e^{-c2^{-m}\Dc}g_m\big (x) \big)^q \,. 
\end{align*}
Therefore, since $\alpha\in(0,1)$, choosing $\eps\in (0,(1-\alpha)/2)$ we obtain
\begin{equation}\label{Sigma-1}
\Sigma_1 (x)
 \lesssim \sum_{m=1}^{+\infty} \sum_{j\ge m/2} 2^{-j(1-\alpha-2\eps)q+mq/2 -\eps mq}
\big(
 e^{-c2^{-m}\Dc}g_m\big (x) \big)^q 
\lesssim \sum_{m=1}^{+\infty} 
\big( 2^{m\alpha/2 } e^{-c2^{-m}\Dc}g_m\big (x) \big)^q \,.
\end{equation}
Analogously, using  H\"older's inequality again, we see that, for
$\eps>0$
\begin{align*}
\bigg( \sum_{m=2j+1}^{+\infty} Mg_{m+1}(x) 
  \bigg)^q  & 
\lesssim \bigg(  \sum_{m=2j+1}^{+\infty}  2^{-\eps mq'} \bigg)^{q/q'} 
\sum_{m=2j+1}^{+\infty}  \big( 2^{\eps m}Mg_{m+1}(x)  \big)^q  \\
& \lesssim 2^{-2\eps jq} 
\sum_{m=2j+1}^{+\infty}  \big( 2^{\eps m}Mg_{m+1}(x)  \big)^q \,,
\end{align*}
so that, if $\eps<\alpha/2$ 
\begin{equation}\label{Sigma-2}
\Sigma_2 (x)
 \lesssim \sum_{m=1}^{+\infty} \sum_{j\le m/2} 2^{\eps mq+jq(\alpha-2\eps)}
\big(
 Mg_{m+1} (x) \big)^q 
\lesssim \sum_{m=1}^{+\infty} 
\big( 2^{m\alpha/2}  Mg_{m+1}(x) \big)^q \,.
\end{equation}
Therefore, from~\eqref{CRTN-argument} we have
\begin{align*}
\big\| \sS^{{\rm loc}, q}_{\alpha}(f-e^{-\Dc}f) \big\|_{L^p(\mu_\chi)} 
& \lesssim \big\| \Sigma_1^{1/q}\big\|_{L^p(\mu_\chi)} + 
\big\| \Sigma_2^{1/q}\big\|_{L^p(\mu_\chi)} \,.
\end{align*}
We first estimate the latter term.  By the Fefferman--Stein 
vector-valued theorem (see \cite{GarciaCuerva-deFrancia} p. 481) with 
the $L^p$-boundedness of the local maximal function, and H\"older's inequality, we
have 
\begin{align*}
\big\| \Sigma_2^{1/q}\big\|_{L^p(\mu_\chi)} 
& = \Big\| \bigg( \sum_{m=1}^{+\infty} 
\big( 2^{m\alpha/2}  Mg_{m+1} \big)^q
\bigg)^{1/q}\Big\|_{L^p(\mu_\chi)}  \\
& \lesssim \Big\| \bigg( \sum_{m=1}^{+\infty} 
\big( 2^{m\alpha/2}  g_{m+1} \big)^q
\bigg)^{1/q}\Big\|_{L^p(\mu_\chi)}  \\
& \lesssim \Big\| \bigg( \sum_{m=1}^{+\infty} 
 2^{m\alpha q/2 -m(q-1)} 
\int_{2^{-m}}^{2^{-m+1}} 
\Big|\frac{\p}{\p t}e^{-t\Dc}f \Big|^q\, dt
\bigg)^{1/q}\Big\|_{L^p(\mu_\chi)} \\
& \lesssim \Big\| \bigg(
\int_0^1 \big( t^{-\alpha/2}| W^{(1)}_t f(x)| \big)^q\, \frac{d t}{t}
\bigg)^{1/q}\Big\|_{L^p(\mu_\chi)} \\
& \lesssim \| f\|_{F^{p,q}_\alpha}  \,. 
\end{align*}

Next, using \eqref{Sigma-1} and applying Proposition 2.4
in~\cite{BPV} and then arguing as before, we estimate
\begin{align*}
\big\| \Sigma_1^{1/q}\big\|_{L^p(\mu_\chi)} 
& \lesssim \Big\| \bigg(
\sum_{m=1}^{+\infty} 
\big( 2^{m\alpha/2 } e^{-c2^{-m}\Dc}g_m  \big)^q
\bigg)^{1/q}\Big\|_{L^p(\mu_\chi)}  
\\
& \lesssim \Big\| \bigg(
\sum_{m=1}^{+\infty} 
\big( 2^{m\alpha/2 } g_m  \big)^q
\bigg)^{1/q}\Big\|_{L^p(\mu_\chi)}  \\
& \lesssim \Big\| \bigg( \sum_{m=1}^{+\infty} 
 2^{m\alpha q/2 -m(q-1)} 
\int_{2^{-m-1}}^{2^{-m}} 
\Big|\frac{\p}{\p t}e^{-t\Dc}f \Big|^q\, dt
\bigg)^{1/q}\Big\|_{L^p(\mu_\chi)} \\
& \lesssim \| f\|_{F^{p,q}_\alpha}  \,. 
\end{align*}
This completes Step 2. \ms

{\em Step 3.} We finish the proof by showing that
$$
\big\|\sS^{{\rm loc}, q}_{\alpha}e^{-\Dc}f
\big\|_{L^p(\mu_\chi)}
\lesssim \| f\|_{L^p(\mu_\chi)}  \,.  \ms
$$

We first notice that for
every 
$x\in  G$ and $y\in B_1$, 
\begin{align*}
|\oD_y(e^{-\Dc}f) (x)|
&\lesssim |y|\sup \big\{ |X_ie^{-\Dc}f(w)|: |w^{-1}x|\leq |y| \,,\
i=1, \dots,\ell\big\}\\
&\lesssim |y|\sup \big\{ |X_ie^{-\Dc}f(w)|: |w^{-1}x|\leq
1\,,\
i=1, \dots,\ell\big\} \,.
\end{align*}
By Lemma~\ref{estimates_pt-chi}  there exists
$t_0>0$ such that for every $w$ such that  $|w^{-1}x|\leq 1$, and
$i=1,\dots,\ell$,  
\begin{align*}
|X_ie^{-\Dc}f(w)|
&=|f * X_i p_1^{\chi}(w)|\\ 
&\leq \int |f(wy^{-1})||X_i p_1^{\chi}(y)|\, d\rho(y)  
\lesssim \int |f(wy^{-1})|  (\delta \chi^{-1})^{1/2}(y) e^{-c|y|^2}  \,
d\rho(y)\\
&\lesssim \int |f(wy^{-1})|   p_{t_0}^{\chi}(y)\, d\rho(y) = \int |f(z)|   p_{t_0}^{\chi}(z^{-1}w) \, d\l(z)\\
&\lesssim \int |f(z)|   p_{t_0}^{\chi}(z^{-1}x) \, d\l(z) \\
& =e^{-t_0\Dc}|f|(x) \,.
\end{align*}
Therefore,
\begin{align*}
\big( \sS^{{\rm loc}, q}_{\alpha}(e^{-\Dc}f(x)) \big)^q
&=\int_0^1\bigg(
  \frac{1}{u^{\alpha}
    V(u)}\int_{|y|<u}|\oD_y(e^{-\Dc}f)(x)|\, d\rho(y)
\bigg)^q\,\frac{d u}{u}\\ 
&\lesssim \int_0^1  \bigg(  \frac{1}{u^\alpha V(u)} 
\int_{|y|<u} u
  e^{-t_0\Dc}|f|(x)  d\rho(y) \bigg)^q\,\frac{d u}{u}\\ 
&\lesssim e^{-t_0\Dc}|f|(x)^q,
\end{align*}
where we used the fact that $\alpha\in (0,1)$. Hence,
$$
\|\sS^{{\rm loc}, q}_{\alpha}e^{-\Dc}f\|_{L^p(\mu_\chi)}\lesssim
\|e^{-t_0\Dc}|f|\|_{L^p(\mu_\chi)}\lesssim
\|f\|_{L^p(\mu_\chi)}\,,   
$$
which completes Step 3, and the proof of the theorem.
\epf
\ms

\subsection{Characterization of Besov norm by differences}

We now prove a characterization of the Besov norm in terms of the
difference operator \eqref{diff-op-eq}.  Its proof  is inspired by
the one of Theorem 1.16 in \cite{Feneuil}, in the case of a
sub-Laplacian without drift on 
a unimodular
group $G$ with  respect to the Haar measure.

We set
\begin{equation}\label{Lpq-functional-eq}
\sA^{p,q}_\alpha(f)= 
\bigg(\int_{|y|\leq 1} \bigg(\frac{\|\oD_y
      f\|_{L^p(\mu_\chi)}}{|y|^\alpha}\bigg) ^q 
\, \frac{d\rho(y)}{V(|y|)}\bigg)^{1/q} \,,  \ms
\end{equation}

\begin{thm}\label{finite-difference-charac-Besov}
Let $\alpha \in (0,1)$ and $p,q\in [1,+\infty]$. Then
\begin{equation}
\|f\|_{B^{p,q}_\alpha}\approx  \sA^{p,q}_\alpha(f) +\|
f\|_{L^p(\mu_\chi)} \,.
\end{equation}
\end{thm}

\proof
We separate the proof in three steps. The first step deals with some
simple integral estimates relying on the classical Schur's test,  
while the second and third steps contain the 
 inequality $\lesssim$ and $\gtrsim $, respectively,  in the statement. 
\ms

{\em Step 1.}
Let $a\in \bbR$, $s\geq 0$, $c>0$ and define the integral kernel
$ K:(0,1)\times G\to[0,+\infty)$ by
$$
K(t,y) =\chi^a(y)  \left( \frac{|y|^2}{t} \right)^s
\frac{V(|y|)}{V(\sqrt{t})} e^{-c|y|^2/t} \,,
$$
and the corresponding integral operator
$$
T_K g(t) = \int_G K(t,y) g(y)\,\frac{d \rho(y)}{V(|y|)} \,.
$$
Then, we show that for all $q\in [1,+\infty)$ 
$$
T_K:  L^q\big(G, d \rho/ V(|\, \cdot\,|)\big) \to L^q \big( (0,1), \, d t/t\big)
$$ 
is bounded. \ms

To this end, it suffices to apply Schur's test, see e.g. \cite[Theorem
6.18]{Folland} after showing that
\begin{equation}
\label{Schur-estimates}
{\rm (a)} \int_G K(t,y)\, \frac{d \rho(y)}{V(|y|)} \lesssim 1\,, 
\qquad  {\rm (b)} \int_0^1 K(t,y) \, \frac{d t}{t} \lesssim 1\,.
\end{equation}

  Observe that
\begin{align*}
\int_G K(t,y)\, \frac{d \rho(y)}{V(|y|)}  = \int_{|y|^2\leq t}
K(t,y)\,\frac{d \rho(y)}{V(|y|)}  + \int_{|y|^2> t}K(t,y)\, \frac{d
  \rho(y)}{V(|y|)}  
=: I + I\!I.
\end{align*}
It is easy to check  that $I\lesssim 1$.  
Moreover, since by \eqref{d-D} for every $j\geq 0$ and $t\in (0,1)$
$$
\frac{V(2^{j+1}\sqrt{t})}{V(\sqrt{t})} 
\lesssim 2^{d(j+1)} e^{D 2^{j+1}}  \,,
$$
we have 
\begin{align*}
I\!I \lesssim 
\sum_{j=0}^\infty \int_{2^j\sqrt{t} \le |y|< {2^{j+1}\sqrt{t}}}
\frac{1}{V(\sqrt{t})} 2^{2js} e^{a c_X  2^{j+1}\sqrt{t}} e^{-c 2^{2j}}\, 
d \rho(y)
\lesssim 1 \,.
\end{align*}
Thus,  condition (a) in~\eqref{Schur-estimates} is
satisfied. 
In order to prove (b), we separate two cases.
If $|y|\geq 1$, then
\begin{align*}
 \int_0^1 K(t,y) \, \frac{d t}{t} \lesssim e^{ac_X |y| -c'|y|^2}V(|y|)
 \int_0^1 \frac{1}{V(\sqrt{t})} e^{-c|y|^2/ t} \,\frac{d
   t}{t}\lesssim 1 
\end{align*}
while, if $|y|\leq 1$, (hence $\chi(y)\lesssim 1$)
\begin{align*}
   \int_0^1 K(t,y) \, \frac{d t}{t} 
& \lesssim \int_{|y|^2}^1 \left(\frac{|y|^2}{t}\right)^{s+d/2}
\, \frac{d t}{t}+ \int_{0}^{|y|^2}
\left(\frac{|y|^2}{t}\right)^{s+d/2} e^{-c |y|^2/t} \,
\frac{d t}{t} \\ 
& \lesssim 1+ \sum_{j=0}^{+\infty}
\int_{2^{-(j+1)}|y|^2}^{2^{-j}|y|^2} 2^{j(s+d/2)}
e^{-c |y|^2/t} \, \frac{d t}{t} \\ 
&   \lesssim 1 \,,
\end{align*}
which proves \eqref{Schur-estimates}.  This completes Step 1.
\ms

{\em Step 2.}
We show that, for $p,q\in [1,+\infty]$ and $\alpha>0$,
\begin{equation}\label{Step2}
\sB^{p,q}_\alpha(f) 
\lesssim \sA^{p,q}_\alpha(f) +
\|f\|_{L^p(\mu_\chi)} \,.\ms
\end{equation}

We claim that there exists $c>0$ such that, for all $f\in L^p(\mu_\chi)$,
\begin{equation}\label{claimteofunc}
\sB^{p,q}_\alpha(f) \lesssim \bigg( \int_G \bigg(
    \frac{\|\oD_y f\|_{L^p(\mu_\chi)} e^{-c|y|^2}}{|y|^\alpha}
  \bigg)^q \, \frac{d \rho(y)}{V(|y|)} \bigg)^{1/q}. 
\end{equation}
Assuming the claim, we prove the estimate~\eqref{Step2}. 
By the claim, it  suffices to prove that
$$
\bigg( \int_G \bigg(
    \frac{\|\oD_y f\|_{L^p(\mu_\chi)} e^{-c|y|^2}}{|y|^\alpha}
  \bigg)^q \, \frac{d \rho(y)}{V(|y|)} \bigg)^{1/q}
 \lesssim
\sA^{p,q}_\alpha(f) + \|f\|_{L^p(\mu_\chi)} \,.
$$
We split the integral on $G$ as $\{|y|<1\} \cup \{|y|\geq 1\}$. On
the one hand, it is easy to see that
$$ 
\bigg( \int_{|y|\leq 1} \bigg( \frac{\|\oD_y f\|_{L^p(\mu_\chi)}
      e^{-c|y|^2}}{|y|^\alpha} \bigg)^q \, \frac{d \rho(y)}{V(|y|)}
\bigg)^{1/q} \leq \sA^{p,q}_\alpha(f) \,,
$$
while since $\|\oD_y f\|_{L^p(\mu_\chi)} \leq  (1 + \chi^{1/p}(y))
\|f\|_{L^p(\mu_\chi)}$, 
\begin{align*} 
\bigg( \int_{|y|\geq 1} \bigg( \frac{\|\oD_y f\|_{L^p(\mu_\chi)}
      e^{-c|y|^2}}{|y|^\alpha} \bigg)^q \, \frac{d \rho(y)}{V(|y|)}
\bigg)^{1/q} & \lesssim  \|f\|_{L^p(\mu_\chi)} \bigg(
  \int_{|y|\geq 1} (1+e^{c(\chi)|y|})e^{-cq|y|^2}\,d \rho(y)
\bigg)^{1/q} \\ 
& \lesssim \|f\|_{L^p(\mu_\chi)} \,,
\end{align*}
since the volume of balls grows at most exponentially (see~\eqref{d-D}).

It remains to prove the claim~\eqref{claimteofunc}. 
 Since $ \int_G \p_t p_t^\chi (y)\, d \rho(y) = \p_t  \int_G p_t^\chi (y)\, d \rho(y) = 0$,
\begin{align*}
 \frac{\p }{\p t}e^{-t\Dc} f(x) & = \int_G f(xy^{-1}) \frac{\p
   p_t^\chi}{\p t}(y)\,d \rho(y) \\ 
&  
= \int_G  \frac{\p p_t^\chi}{\p t}(y) \oD_y f(x)\, d\rho(y) \,.
\end{align*}
Thus,
\begin{equation}\label{dtptchi}
\Big\| \frac{\p }{\p t} e^{-t\Dc}f\Big\|_{L^p(\mu_\chi)}
\leq \int_G \Big|\frac{\p p_t^\chi}{\p t}(y)
\Big| \|\oD_yf\|_{L^p(\mu_\chi)}\,d \rho(y) \,.
\end{equation}
By Lemma~\ref{estimates_pt-chi} (v),~\eqref{dtptchi} and   
\cite[Lemma 3.3]{BPV} 
\begin{align*}
   \sB^{p,q}_\alpha(f)^q 
&   \approx \int_0^1 \bigg(
  t^{1-\alpha/2} 
\Big\|
\Dc e^{-\frac{t}{2}\Dc}
\frac{\p }{\p t} e^{-\frac{t}{2}\Dc} f\Big\|_{L^p(\mu_\chi)}
\bigg)^q \, \frac{d t}{t}\\
   & \lesssim \int_0^1 \bigg( t^{1- \alpha/2} \int_G
     (\delta \chi^{-1})^{1/2}(y) t^{-(d+2)/2} e^{-b |y|^2/t}  \|\oD_y
     f\|_{L^p(\mu_\chi)}\,d \rho(y) \bigg)^q \, \frac{d t}{t} \\ 
& \lesssim \int_0^1 \bigg( t^{-\alpha/2} \int_G  
 (\delta \chi^{-1})^{1/2}(y) 
  t^{-d/2} e^{-b' |y|^2/t}  \|\oD_y f\|_{L^p(\mu_\chi)}
  e^{-b'|y|^2}d \rho(y) \bigg)^q \, \frac{d t}{t} \\ 
& = \int_0^1 \bigg( \int_G K(t,y)  g(y) \, \frac{d \rho(y)}{V(|y|)}
\bigg)^q \, \frac{d t}{t} 
\end{align*}
with 
$$
b' = \frac{b}{2}, \quad g(y) = \frac{\|\oD_y
  f\|_{L^p(\mu_\chi)}}{|y|^\alpha} e^{-b'|y|^2}, 
\quad K(t,y) = \frac{V(|y|)}{t^{d/2}} \bigg(\frac{|y|^2}{t} 
\bigg)^{\alpha/2}  (\delta \chi^{-1})^{1/2}(y)  e^{-b' |y|^2/t} \,.
$$
By Step 1 we obtain 
\begin{align*}
   \sB^{p,q}_\alpha(f)^q 
& \lesssim  \int_G |g(y)|^q \, \frac{d \rho(y)}{V(|y|)}  = \int_G
\bigg( \frac{\|\oD_y f\|_{L^p(\mu_\chi)} e^{-b'|y|^2}}{|y|^\alpha}
\bigg)^q \, \frac{d \rho(y)}{V(|y|)}. 
\end{align*}
The claim is proved, and Step 2 is complete. \ms

{\em Step 3.} We prove that, for
$p,q\in [1,+\infty]$ and $\alpha \in (0,1)$, 
$$
\sA^{p,q}_\alpha(f) \lesssim \sB^{p,q}_\alpha(f) +
\|f\|_{L^p(\mu_\chi)} \,. \ms
$$
 
We write again
 $f$ as $f = (f-e^{-\Dc} f) + e^{-\Dc} f$ and
 decompose the $(f-e^{-\Dc} f)= \sum_{m=1}^{\infty} f_m$ as in
 \eqref{m-decomposition}. Then, using also  \eqref{m-decomposition-2}, we have 
\begin{align*}
f_m 
& = 
\int_{2^{-m}}^{2^{-m+1}}
\frac{\p}{\p t}e^{-t\Dc}f \, dt =
\int_{2^{-m}}^{2^{-m+1}} \Dc e^{-t\Dc} f\, dt 
= 2\int_{2^{-m-1}}^{2^{-m}} \Dc e^{-2t\Dc}
   f\,d t \\
& = 2 e^{-2^{-m-1}\Dc}\int_{2^{-m-1}}^{2^{-m}}
   e^{-(t-2^{-m-1})\Dc}\Dc e^{-t\Dc} f\, dt
   =:2 e^{-2^{-m-1}\Dc}  h_m \,.
\end{align*}

We set
$$
\sigma_m 
= \int_{2^{-m-1}}^{2^{-m}} 
\Big\|\frac{\p }{\p t}e^{-t\Dc} f\Big\|_{L^p(\mu_\chi)}\, dt 
$$
and observe that 
\begin{equation}\label{fncn1}
\| f_m\|_{L^p(\mu_\chi)} \lesssim \sigma_{m+1} \,.
\end{equation}
Hence,
\begin{equation}\label{fncn2}
\| \oD_y f_m\|_{L^p(\mu_\chi)} \lesssim \big(1+\chi^{1/p}(y)\big)
\|f_m\|_{L^p(\mu_\chi)} 
\lesssim \sigma_{m+1} \,.
\end{equation}

By \cite[Lemma 3.3]{BPV} it follows that for $i=1,\dots,\ell$, 
\begin{align*}
   \|X_i f_m\|_{L^p(\mu_\chi)} 
& \lesssim 2^{m/2} \| h_m\|_{L^p(\mu_\chi)}  
\lesssim 2^{m/2} \int_{2^{-m-1}}^{2^{-m}} 
\Big\|e^{-(t-2^{-m-1})\Dc}\frac{\p }{\p t}e^{-t\Dc} f
\Big\|_{L^p(\mu_\chi)} \, dt  
\lesssim 2^{m/2} \sigma_m \,.
\end{align*}
Therefore,
\begin{equation}\label{eqdist}
  \|\oD_y f_m\|_{L^p(\mu_\chi)}  
\le |y| \sum_{i=1}^\ell \|X_i f_m\|_{L^p(\mu_\chi)}  
\lesssim |y|2^{m/2} \sigma_m \,.
\end{equation}
Since $\sigma_m \le 2\sigma_{m+1}$, \eqref{eqdist},~\eqref{fncn1} and~~\eqref{fncn2} imply 
$$
\|\oD_y f_m\|_{L^p(\mu_\chi)} 
\lesssim 
\begin{cases}
|y|2^{m/2} \sigma_m & \text{ if } |y|^2 < 2^{-m} \\ 
\sigma_{m+1} & \text{ if } |y|^2 \ge 2^{-m}
\end{cases} \,.
$$
Therefore, by  Lemma 4.3 (i) in~\cite{BPV}, we have
\begin{align*}
  \Big( \sA^{p,q}_\alpha(f- e^{-\Dc} f) \Big)^q 
& \lesssim \sum_{j=1}^{\infty} \int_{2^{-j}\leq |y|^2< 2^{-j+1}}
2^{jq\alpha/2} \Big(\sum_{m=1}^{\infty} 
\big\|\oD_y f_m \big\|_{L^p(\mu_\chi)} \Big) ^q \, \frac{d \rho(y)}{V(|y|)} \\ 
& \lesssim \sum_{j=1}^{\infty} 2^{jq\alpha/2}
\bigg(\sum_{m=1}^j 2^{(m-j)/2} \sigma_m 
+ \sum_{m=j+1}^{+\infty} \sigma_{m+1} \bigg) ^q  \\ 
& \lesssim \sum_{j=1}^{\infty} 2^{jq(\alpha-1)/2}
\bigg(\sum_{m=1}^{\infty} 2^{\min \{j,m\}/2} \big(
\sigma_{m+1} + \sigma_m\big) \bigg) ^q  \\
& \lesssim \sum_{m=1}^{\infty} \Big( 2^{m\alpha/2}
\big(\sigma_{m+1} + \sigma_m\big) \Big) ^q  \\
& \lesssim 
\sum_{m=0}^{\infty} \bigg( 2^{m\alpha/2} \int_{2^{-m-1}}^{2^{-m}}
\Big\| \frac{\p }{\p t} e^{-t\Dc} f\Big\|_{L^p(\mu_\chi)} \,
dt \bigg)^q \\
& \lesssim \sum_{m=0}^{\infty} 2^{mq\alpha/2}  2^{-m(q-1)}
\int_{2^{-m-1}}^{2^{-m}} \Big\|\frac{\p }{\p t}e^{-t\Dc} f
\Big\|_{L^p(\mu_\chi)}^q\, dt  \\
& \lesssim \sum_{m=0}^{\infty}  \int_{2^{-m-1}}^{2^{-m}}\Big(
  t^{1-\alpha/2} \Big\|\frac{\p }{\p
      t}e^{-t\Dc} f \Big\|_{L^p(\mu_\chi)}\Big)^q \, \frac{d
      t}{t} =  \big( {\mathscr B}^{p,q}_{\alpha}(f)\big)^q \,.
\end{align*}
Therefore, by~\eqref{g-fct-est} we have 
$$
\sA^{p,q}_\alpha(f- e^{-\Dc} f) 
\lesssim
\sB^{p,q}_\alpha(f) \,.
$$

It remains to estimate $\sA^{p,q}_\alpha(e^{-\Dc} f)  $. As in~\eqref{eqdist}, one can see that
$\|\oD_y e^{-\Dc} f\|_{L^p(\mu_\chi)} \lesssim |y|\|f\|_{L^p(\mu_\chi)}$
so that, using the decomposition of the integral as sum of
integrals over annuli, 
\begin{align*}
    \sA^{p,q}_\alpha(e^{-\Dc} f)  
& \leq \|f\|_{L^p(\mu_\chi)} \bigg( \int_{|y|\leq 1} |y|^{q(1-\alpha)}
\, \frac{d \rho(y)}{V(|y|)} \bigg)^{1/q}  \lesssim
\|f\|_{L^p(\mu_\chi)} \,.
\end{align*}
The proof of Step 3 is complete.
This concludes the proof of Theorem~\ref{finite-difference-charac-Besov}.
\epf

\ms

\section{A density result}\label{density-sec}
\ms

The main goal of this section is to show that the smooth functions
with compact support are dense in the Triebel--Lizorkin and Besov
spaces on $G$. This is the analogue of the classical density result in the Euclidean setting \cite{Triebel}; we refer the reader to \cite{HMY, HKT, Triebel-II} for its counterpart in other settings.

To prove our density results we shall use the following version of
Young's inequality. 
\begin{lem}\label{Young-ineq-lem}
If $\eta$ has support in $B_1$, then
\begin{equation}\label{Young-ineq}
\| \eta*f\|_{L^p(\mu_\chi)} \lesssim 
\| \eta\|_{L^1(\l)} \| f\|_{L^p(\mu_\chi)} \,.
\end{equation}
\end{lem}

\proof
Suppose $\eta$ has compact support
contained in $B_1$, then by Minkowski's integral inequality,
\begin{align*}
\|\eta* f \|_{L^p(\mu_\chi)}
&= \Big\| \int_G \eta(y^{-1}) f(y\cdot)\, d\rho(y) \Big\|_{L^p(\mu_\chi)}\le
\int_{B_1} |\eta(y^{-1})| \|  f(y\cdot)\|_{L^p(\mu_\chi)}\, d\rho(y) \,.
\end{align*}
Now,
\begin{align*}
\|  f(y\cdot)\|_{L^p(\mu_\chi)}^p
& = \int_G |f(yx)|^p (\chi\delta^{-1})(x)\, d\l(x) =  
(\chi^{-1} \delta)(y) \int_G |f(x)|^p (\chi\delta^{-1})(x)\, d\l(x) \\
& =  (\chi^{-1}\delta)(y) \int_G |f(x)|^p \, d\mu_\chi(x) \,.
\end{align*}
Therefore,
\begin{align*}
\|\eta* f \|_{L^p(\mu_\chi)}  
& \le \Big( \int_{B_1} |\eta(y^{-1})| \big(  (\chi^{-1}
\delta)a(y) \big)^{1/p} \, d\rho(y) \Big) \|f \|_{L^p(\mu_\chi)} \\
&  =  
\Big( \int_{B_1} |\eta(y)| \big( (\chi\delta^{-1})(y) \big)^{1/p} \, d\l(y) \Big) \|f \|_{L^p(\mu_\chi)}\\
& \le \sup_{B_1} \big( \chi\delta^{-1} \big)^{1/p} \| \eta\|_{L^1(\l)} \| f\|_{L^p(\mu_\chi)} \,.
\end{align*}
The conclusion follows. \qed

\ms

We now prove the density result.  

\begin{thm}\label{density-lem}
Let $\alpha>0$, $p,q\in(1,\infty)$ and let
$X^{p,q}_\alpha$ denote either space $F^{p,q}_\alpha$ or
$B^{p,q}_\alpha$. Then, $C_c^\infty(G)$ is dense in $X^{p,q}_\alpha$.
\end{thm}

\proof
We begin by observing that Lemma \ref{pt-Schwartz-cont} implies that
$\cS$, hence $C^\infty_c$, are contained in $X^{p,q}_\alpha$.  Indeed, if
$m,n\in\bbN$ with $m>[\alpha/2]$ and $n$ to be chosen, for 
$\vp \in \cS$, using \eqref{Sch-est} we have 
$$
|W_t^{(m)} \vp(x)| \le t^m e^{-n|x|} \sum_{|J|\le 2m} \cN_{J,n} (\vp) \,.
$$
Then, in order to estimate the norms in~\eqref{TL-norm}
and~\eqref{Besov-norm}, 
it suffices to chose $n$ large enough so that $\int_G e^{-pn|x|}
d\mu_\chi(x)$ is finite.

{\em Step 1.}  We first prove 
that we can approximate functions in $F^{p,q}_\alpha$ with functions
having compact support when $\alpha\in(0,1)$.   

Let $\eta\in C^\infty_c(B_1)$, $\eta=1$ on
$B_{1/2}$.  Given any $R>2$, define $\eta_R= {\mathbf
  1}_{B_R}*\eta$. Then, 
$\eta_R\in C^\infty_c(B_{R+1})$ and $\eta_R(x)=1$ on
$B_{R-1}$.   We 
 observe that $\|X_J\eta_R\|_\infty \lesssim 1$, for $|J|\le n$.
In order to show such bounds, we apply Young's inequality~\eqref{Young-ineq}. 
For any $J\in \sJ^k$, with $k\le n$, we have
$$
\|X_J \eta_R\|_\infty = \|  {\mathbf
  1}_{B_R}*X_J\eta \|_\infty \le \|X_J \eta \|_{L^1(\l)}\lesssim 1\,.
$$

Then, let $f\in F^{p,q}_\alpha$ be given.  
We shall estimate the norm
$\| f - f\eta_R\|_ {F^{p,q}_\alpha}$ using the $\sS^{{\rm
    loc},q}_\alpha$-functional. Since $\alpha\in(0,1)$, $f\in L^p(\mu_\chi)$ and $\| f -
f\eta_R\|_{L^p(\mu_\chi)} \to 0$ as $R\to+\infty$. 
Next, we show that also
\begin{align}
\big\| \sS^{{\rm loc},q}_\alpha 
\big(  f - f\eta_R\big) \big\|_{L^p(\mu_\chi)} \to 0 \,, \label{appr-comp}
\end{align}
as $R\to+\infty$. 
Set $\z_R=1-\eta_R$ and observe that $\z_R$ vanishes
identically
on $B_{R-1}$.  
Then, we have
\begin{align*}
\big[ \sS^{{\rm loc},q}_\alpha 
\big(  f - f\eta_R\big)(x) \big]^q
&  = \int_0^1 \Big[ \frac{1}{u^\alpha V(u)} \int_{|y|\le u}
|\z_R(xy^{-1}) f(xy^{-1}) - \z_R(x)f(x)|\, d\rho(y) \Big]^q\, \frac{du}{u}
\\
& \lesssim \int_0^1 \Big[ \frac{1}{u^\alpha V(u)} \int_{|y|\le u}
\z_R(xy^{-1}) |f(xy^{-1}) - f(x)|\, d\rho(y) \Big]^q\, \frac{du}{u}\\
& \qquad
+  |f(x)|^q\int_0^1 \Big[ \frac{1}{u^\alpha V(u)} \int_{|y|\le u}
|\z_R(xy^{-1}) - \z_R(x)|\, d\rho(y) \Big]^q\, \frac{du}{u} \\
& \lesssim \int_0^1 \Big[ \frac{1}{u^\alpha V(u)} \int_{|y|\le u}
\z_R(xy^{-1}) |f(xy^{-1}) - f(x)|\, d\rho(y) \Big]^q\, \frac{du}{u}\\
& \qquad
+  |f(x)|^q \int_0^1 \Big[ \frac{1}{u^\alpha V(u)} \int_{|y|\le u}
|y|\sup_{z\in B(x,1)} |\nabla \z_R(z)|\, d\rho(y) \Big]^q\, \frac{du}{u} \\
& \lesssim {\mathbf 1}_{\{|x|\ge R-2\}}(x) \Big( \big[\sS^{{\rm loc},q}_\alpha
f(x)\big]^q + |f(x)|^q  \Big)\,,
\end{align*}
since $\z_R(z)=0$ if $|z|\le R-1$, and $\|\z_R\|_\infty=1$.  Therefore,
\begin{align*}
\big\| \sS^{{\rm loc},q}_\alpha 
\big(  f - f\eta_R\big) \big\|_{L^p(\mu_\chi)}^p
& \lesssim \int_{\{|x|\ge R-2\}} \big[\sS^{{\rm loc},q}_\alpha
f(x)\big]^p\, d\mu_\chi(x)  + \int_{\{|x|\ge R-2\}} |f(x)|^p\, d\mu_\chi(x)\,.
\end{align*}
This proves \eqref{appr-comp} and therefore we can approximate any element of 
$F^{p,q}_\alpha$ with elements with compact support.\ms

{\em Step 2.}    Using the charaterization of the norm in
$B^{p,q}_\alpha$ by finite
differences for $\alpha\in(0,1)$, Theorem~\ref{finite-difference-charac-Besov}, we prove
that we can approximate functions in $B^{p.q}_\alpha$ with functions
with compact support.  

Let $f \in B^{p.q}_\alpha$ be given and let
$\eta_R$ and $\z_R$ be as in Step 1.  Then, for $y\in B_1$ and $x\in
G$ we write
\begin{align*}
\oD_y(f\z_R) (x) & = (f\z_R)(xy^{-1}) - (f\z_R)(x) \\
& = \z_R(xy^{-1}) \big[ f(xy^{-1}) - f(x) \big] + f(x) \big[
  \z_R(xy^{-1}) -\z_R(x) \big] =: F_R(x,y) + G_R(x,y) \,. 
\end{align*}
We observe that 
\begin{align*}
\bigg(\int_{|y|\leq 1} \bigg(
\frac{\|G_R(\cdot,y)\|_{L^p(\mu_\chi)}}{|y|^\alpha} \bigg) ^q 
\, \frac{d\rho(y)}{V(|y|)}\bigg)^{1/q}
& \lesssim \bigg(\int_{|y|\leq 1} 
|y|^{(1-\alpha)q}
\, \frac{d\rho(y)}{V(|y|)}\bigg)^{1/q}  
\|{\mathbf 1}_{\{|x|\ge  R-2\}}f \|_{L^p(\mu_\chi)} \\
& \lesssim \|{\mathbf 1}_{\{|x|\ge R-2\}}f \|_{L^p(\mu_\chi)} \,,
\end{align*}
that tends to 0 as $R\to+\infty$.  Next, we observe that
$|F_R(x,y)| \le {\mathbf 1}_{\{|x|\ge R-2\}} |\oD_y f(x)|\le |\oD_y f(x)|$, so that 
\begin{itemize}
\itb $\displaystyle{ 
\|F_R(\cdot,y)\|_{L^p(\mu_\chi)}
\le \Big(\int_{\{|x|\ge R-2\}} |\oD_y f(x)|^p \, d\mu_\chi(x) \Big)^{1/p}
\to 0}$,  as $R\to+\infty$;\smallskip 
\itb
$\|F_R(\cdot,y)\|_{L^p(\mu_\chi)} \le \|\oD_y
f\|_{L^p(\mu_\chi)}$, which is indipendent of $R$.
\end{itemize}
Lebesgue's theorem now gives that 
$$
\bigg(\int_{|y|\leq 1} \bigg(
\frac{\|F_R(\cdot,y)\|_{L^p(\mu_\chi)}}{|y|^\alpha} \bigg) ^q 
\, \frac{d\rho(y)}{V(|y|)}\bigg)^{1/q}
\to 0
$$
as $R\to+\infty$.
Hence, recalling~\eqref{Lpq-functional-eq}, we have 
\begin{align*}
\sA^{p,q}_\alpha(f-f\eta_R)
& \le  
\bigg(\int_{|y|\leq 1} \bigg(\frac{\|F_R(\cdot,y)\|_{L^p(\mu_\chi)}}{|y|^\alpha}\bigg) ^q 
\, \frac{d\rho(y)}{V(|y|)}\bigg)^{1/q} \\
& \qquad \qquad + \bigg(\int_{|y|\leq 1} \bigg(\frac{\|G_R(\cdot,y)\|_{L^p(\mu_\chi)}}{|y|^\alpha}\bigg) ^q 
\, \frac{d\rho(y)}{V(|y|)}\bigg)^{1/q} \to 0 \,,
\end{align*}
as $R\to+\infty$.  This completes Step 2.
\ms

{\em Step 3.}  
We select a smooth
approximation of the identity.  Precisely, for $0<\kappa\le1/2$, select
$\eta_\kappa\in
C^\infty_c(B_\kappa)$, $\eta_\kappa\ge0$, and
$\|\eta_\kappa\|_{L^1(\l)}=1$.  Moreover, we require that
$\eta_\kappa\lesssim V(2\kappa)^{-1}$ (where the constant 
does not depend on $\kappa$).\footnote{This can achieved as follows.
  Let $\tilde\eta\in C^\infty_c$, $0\le\eta\le1$, $\supp\eta\subseteq
  B_{\kappa/2}$, $\tilde\eta=1$ on $B_{\kappa/4}$, and define
  $\eta_\kappa= C_\kappa V(2\kappa)^{-1} {\mathbf
    1}_{B_{\kappa/4}}*\tilde\eta$.  By requiring that
  $\|\eta_\kappa\|_{L^1(\l)}=1$, we obtain that $C_\kappa\approx
  C^{-1}$, where $C$ is the local doubling constant.}
We then have
\begin{equation}\label{approx-ident-max-fnct-control}
|\eta_\kappa*f(x)|\lesssim Mf(x) \,,
\end{equation}
where $M$ is the local maximal function, defined in~\eqref{MR}.
Indeed, we estimate
$$
|\eta_\kappa*f(x)| =
\big| \int_G \eta_\kappa (xy^{-1}) f(y)\, d\rho(y) \big|
\lesssim \frac{1}{V(2\kappa)}  \int_{B(x,\kappa)} |f(y)|\, d\rho
\lesssim Mf(x)\,,
$$
as claimed.

Arguing as  in \cite{Folland-CAHA} Proposition 2.44, by
Minkowski's integral inequality, for any $g\in L^p(\mu_\chi)$,
\begin{align*}
\| \eta_\kappa*g - g\|_{L^p(\mu_\chi)} 
& = 
\big\| \int_G \eta_\kappa(y) (\tau_y g - g)\, d\l(y)
\big\|_{L^p(\mu_\chi)} 
\le  \int_G \eta_\kappa(y) \|\tau_y g - g \|_{L^p(\mu_\chi)} \, d\l(y)
\\
& \le \sup_{|y|\le\kappa} \|\tau_y g - g \|_{L^p(\mu_\chi)}\,,
\end{align*}
which tends to 0, as $\kappa\to0$. 
Next, notice that, by left invariance
$
W^{(m)}_t (\eta_\kappa*f) 
=\eta_\kappa* W^{(m)}_tf
$.
Moreover, if $f\in\cS'$ has compact support, then $\eta_\kappa*f\in C^\infty_c(G)$.\ms 

{\em Step 4.} We complete the proof that $C_c^\infty$ is dense in the
Triebel--Lizorkin spaces, $F^{p,q}_\alpha$ in the case
$\alpha\in(0,1)$. 
 To this end, let $f\in F^{p,q}_\alpha$ have compact support so that
$\eta_\kappa*f \in C_c^\infty$.

By Theorem~\ref{teo-equiv2} we have
\begin{align*}
\| \eta_\kappa*f - f\|_{F^{p,q}_\alpha} 
& \lesssim  \Big \| \Big( \sum_{j=0}^{+\infty} \big(2^{j\alpha/2} 
|\eta_\kappa*W^{(m)}_{2^{-j}} f -W^{(m)}_{2^{-j}} f|\big)^q
\Big)^{1/q} \Big\|_{L^p(\mu_\chi)}  
+ \|\eta_\kappa* f -  f\|_{L^p(\mu_\chi)}  \,. 
\end{align*}
We only need to estimate the first term in the right hand side above.
Observe that, since $f$ has compact support, $W^{(m)}_tf\in \cS$, for
all $t\in(0,1)$.  Hence, \cite{Folland-CAHA} Proposition 2.44 gives
that, for each $j$ fixed
$$
\|
\eta_\kappa*W^{(m)}_{2^{-j}} f -
W^{(m)}_{2^{-j}} f \|_\infty \to 0
\qquad\text{as}\quad \kappa\to 0\,.
$$
We wish to apply Lebesgue's theorem to the inner sum. We observe that 
by \eqref{approx-ident-max-fnct-control} we have that 
$$
\big( 2^{j\alpha/2} 
|\eta_\kappa*W^{(m)}_{2^{-j}} f -W^{(m)}_{2^{-j}} f|\big)^q
\lesssim
\big( 2^{j\alpha/2} \big[ M (W^{(m)}_{2^{-j}} f) + |W^{(m)}_{2^{-j}} f|
\big]  \big)^q\,,
$$
which is (independent of $\kappa$ and) summable by
\cite{GarciaCuerva-deFrancia} p. 481.  Thus, the inner sum tends to 0,
as $\kappa\to0$, for every $x\in G$, that is, the family of
vector-valued 
functions $S_\kappa$, 
$$
S_\kappa:= \big(2^{j\alpha/2} \eta_\kappa*W^{(m)}_{2^{-j}} f \big)_j :G\to \ell^q
$$
as $\kappa\to0$ converges pointwise to the vector-valued function
$S:G\to\ell^q$, where
$S:=\big(2^{j\alpha/2} W^{(m)}_{2^{-j}} f\big)_j$. 
Since, as before, for $0<\kappa<1$, 
\begin{align*} 
\| S_\kappa(x)\|_{\ell^q} &
\lesssim 
\Big(\sum_{j=0}^{+\infty} \big(2^{j\alpha/2} M
(W^{(m)}_{2^{-j}} f(x) \big)^q\Big)^{1/q} \in L^p\big(
(G,\ell^q),\mu_\chi\big)\,. 
\end{align*}
We can apply Lebesgue's theorem to obtain that $S_\kappa \to S$ in 
$L^p\big(
(G,\ell^q),\mu_\chi\big)$, that is, 
$$
 \Big \| \Big( \sum_{j=0}^{+\infty} \big(2^{j\alpha/2} 
|\eta_\kappa*W^{(m)}_{2^{-j}} f -W^{(m)}_{2^{-j}} f|\big)^q
\Big)^{1/q} \Big\|_{L^p(\mu_\chi)}  \to 0
\qquad\text{as}\quad \kappa\to 0\,,
$$
as we wished to show.  Hence,  $C_c^\infty$ is dense in
$F^{p,q}_\alpha$, for $p,q\in(1,\infty)$, $\alpha\in(0,1)$.\ms

{\em Step 5.} 
We complete the proof that $C_c^\infty$ is dense in
$B^{p,q}_\alpha$ in the case $\alpha\in(0,1)$. Let $f\in B^{p,q}_\alpha$ have compact support so that
$\eta_\kappa*f \in C_c^\infty$.  We then have
\begin{align*}
\| \eta_\kappa*f - f\|_{B^{p,q}_\alpha} 
& \lesssim  \Big( \int_0^1 \big(t^{-\alpha/2} \| \eta_\kappa*W^{(m)}_tf -
W^{(m)}_tf\|_{L^p(\mu_\chi)} \big)^q \, \frac{dt}{t} \Big)^{1/q}  + \|
\eta_\kappa* f -  f\|_{L^p(\mu_\chi)} \,.
\end{align*}
Now, $\| \eta_\kappa* f - 
f\|_{L^p(\mu_\chi)}$ and 
 $\| \eta_\kappa*W^{(m)}_tf - W^{(m)}_tf\|_{L^p(\mu_\chi)} \to 0$, as
 $\kappa\to0$, the latter term for each $t\in(0,1)$ fixed.
Using Young's inequality~\eqref{Young-ineq} we see that
$$
\| \eta_\kappa*W^{(m)}_tf - W^{(m)}_tf\|_{L^p(\mu_\chi)} 
\lesssim \|W^{(m)}_tf\|_{L^p(\mu_\chi)} \,,
$$
so that 
we may use
  Lebesgue's theorem to obtain that
$$
\Big( \int_0^1 \big(t^{-\alpha/2} \| \eta_\kappa*W^{(m)}_tf -
W^{(m)}_tf\|_{L^p(\mu_\chi)} \big)^q \, \frac{dt}{t} \Big)^{1/q} 
\to 0\qquad\qquad \text{as} \quad \kappa\to0\,.
$$

This gives that 
$\| \eta_\kappa*f - f\|_{B^{p,q}_\alpha} \to 0$, as $\kappa\to0$,  
and completes the proof that $C_c^\infty$ is dense in
$B^{p,q}_\alpha$, for $p,q\in(1,\infty)$, $\alpha\in(0,1)$.
\ms

{\em Step 6.} We now prove that $C_c^\infty$ is dense in
$X^{p,q}_{\alpha+1}$, 
for $p,q\in(1,\infty)$, $\alpha\in(0,1)$.  We use the 
recursive characterizations of the spaces $X^{p,q}_{\alpha+1}$, see
Theorem 4.5 in~\cite{BPV}, that is, $f\in
X^{p,q}_{\alpha+1}$ if and only if $f\in X^{p,q}_\alpha$ and $X_j f\in 
X^{p,q}_\alpha$.  
Let $\eps>0$. 
By the arguments in Steps 1 and 2, using the same notation, 
 there exists $R>0$ sufficiently large such that 
$$
\| f-f\eta_R\|_{X^{p,q}_\alpha} <\eps\,,\qquad \text{and} \quad
\| X_jf-(X_jf)\eta_R\|_{X^{p,q}_\alpha} <\eps/2 \,,
$$
for $j=1,\dots,\ell$. 
Since $X_j\eta_R$ is a $C^\infty_c$ function vanishing in the ball
$B_R$, by the arguments involving $\z_R$ in Steps 1 and 2 we can also assume that
$\|fX_j\eta_R\|_{X^{p,q}_\alpha} <\eps/2$, $j=1,\dots,\ell$.
Therefore,
$$
\| f-f\eta_R\|_{X^{p,q}_\alpha} <\eps\,,\qquad \text{and} \quad
\| X_jf-X_j(f\eta_R)\|_{X^{p,q}_\alpha} <\eps \,,
$$
that is, $\| f-f\eta_R\|_{X^{p,q}_{\alpha+1}} \lesssim \eps$.  

Next, given $f\in X^{p,q}_{\alpha+1}$ having compact support, let
$\{\eta_\kappa\}$, $0<\kappa<1$ be the approximation of the identity
of Step 3, and consider $\eta_\kappa*f$.  Then, by Steps 4 and 5,
$\eta_\kappa*f \to f$, $X_j(\eta_R*f) =\eta_\kappa*X_jf \to X_j f$, $j=1,\dots,\ell$,
in $X^{p,q}_\alpha$, as $\kappa\to0$.  This implies that $\eta_R f
\to f$ in $X^{p,q}_{\alpha+1}$, as $\kappa\to0$. This shows that
$C_c^\infty$ is dense in $ X^{p,q}_{\alpha+1}$, for
$\alpha\in(0,1)$. \ms

{\em Step 7.} We now finish the proof.  Arguing as in the previous
step, we obtain that  $C_c^\infty$ is dense in $ X^{p,q}_\alpha$ for
all $\alpha\in\bbR\setminus\bbN$.  Let $n$ be a positive integer, and 
$\theta\in(0,1)$. 
By the complex interpolation
results of Theorem 6.1 in~\cite{BPV}, $X^{p,q}_n
= ( X^{p,q}_{n-\theta}, X^{p,q}_{n+\theta})_{[1/2]}$, with
$\theta\in(0,1)$.
  By Theorem
4.2.2 in \cite{Bergh-Lofstrom} $X^{p,q}_{n-\theta} \cap
X^{p,q}_{n+\theta}$ is dense in $X^{p,q}_n$.  Let $f\in 
X^{p,q}_{n-\theta} \cap
X^{p,q}_{n+\theta}$ and let $\eta\in C^\infty_c$ such that
$\| f-\eta\|_{X^{p,q}_{n+\theta}}<\eps$, so that also 
$\| f-\eta\|_{X^{p,q}_{n-\theta}}<\eps$.  By the interpolation
inequality  (see e.g. ~\cite[Theorem 4.7.1, p.\
102]{Bergh-Lofstrom}  
and~\cite[p.\ 49]{Bergh-Lofstrom}) 
$$
\| f-\eta\|_{X^{p,q}_n} \le \| f-\eta\|_{X^{p,q}_{n-\theta}}^{1/2}
\| f-\eta\|_{X^{p,q}_{n+1-\theta}}^{1/2} < \eps
$$
we obtain that $C^\infty_c$ is also dense in $X^{p,q}_n$.
\epf

\ms

\section{Isomorphisms of Triebel--Lizorkin and Besov
  spaces}\label{iso-sec} 
\ms 

Goal of this section is to prove that Bessel potentials provide
isomorphisms in both the Triebel--Lizorkin and Besov scales and that a simplified version of local Riesz transforms 
is bounded on  both the Triebel--Lizorkin and Besov spaces.
We continue to denote by $X^{p,q}_\alpha$ either space
$F^{p,q}_\alpha$, or $B^{p,q}_\alpha$.
We begin by showing that for all $c\ge0$, the fractional powers of $\Dc+cI$ are
bounded on the spaces $X^{p,q}_\alpha$.  Precisely, we prove the following.

\begin{lem}\label{frac-sublap-bdd-lem}
Let $p,q\in (1,+\infty)$, $\alpha\ge0$ and let $\gamma>0$.  Then, for
all $c\ge0$, 
$$
(\Dc+cI)^{\gamma/2}: X^{p,q}_{\alpha+\gamma}\to X^{p,q}_\alpha
$$
is bounded.
\end{lem}

\proof
If $\beta>0$, then for $\tau>0$ we have
$$
\tau^{-\beta}
=\frac{1}{\Gamma(\beta)} \int_0^{+\infty} s^{\beta-1}e^{-\tau s}\, ds
\,,
$$
so that, if $\gamma>0$ and $k$ is an integer, $k\ge[\gamma/2]+1$,
\begin{align}
(\Dc+cI)^{\gamma/2} & = (\Dc+cI)^{-(k-\gamma/2)} (\Dc+cI)^k \notag \\
& =
 \frac{1}{\Gamma(k-\gamma/2)} \int_0^{+\infty}
s^{k-\gamma/2}e^{-s\Dc}(\Dc+cI)^k \, \frac{ds}{s} 
\,. \ms \label{subord-1}
\end{align}

We first prove the result for $c=0$.  The case $c>0$ will then follow
easily from the former case.    \ms

{\em Step 1.} We prove that for all $f\in X^{p,q}_\alpha$
\begin{equation}\label{Besov-TL-ineq-1}
\big\| e^{-\frac12\Dc} \Dc^{\gamma/2} f
\big\|_{L^p(\mu_\chi)} \lesssim \| e^{-\frac14 \Dc}
f\|_{L^p(\mu_\chi)} \,. \ms
\end{equation}
Using \eqref{subord-1} (with $c=0$), Lemma 3.3 in~\cite{BPV}, the Cauchy--Schwarz
 inequality and the boundedness of the $g$-function~\eqref{g-fct-est},  
we notice that 
\begin{align*}
\big\| e^{-\frac12\Dc} \Dc^{\gamma/2} f
\big\|_{L^p(\mu_\chi)} 
& \lesssim \int_0^{3/4} s^{k-\alpha/2} \big\| e^{-(s+1/4)\Dc} \Dc^k
e^{-\frac 14 \Dc}  f\big\|_{L^p(\mu_\chi)} \, \frac{ds}{s} \\
& \qquad\qquad 
+ \Big\| \bigg( \int_{3/4}^{+\infty} \big| (s\Dc)^k e^{-s\Dc} e^{-\frac 14 \Dc}
f\big|^2 \, \frac{ds}{s} \bigg)^{1/2} \Big\|_{L^p(\mu_\chi)} 
\\ 
& \lesssim \int_0^1\frac{s^{k-\alpha/2}}{(s+1/4)^k} \big\| e^{-\frac 14 \Dc}
f\big\|_{L^p(\mu_\chi)} \, \frac{ds}{s} + \big\| g_k(e^{-\frac 14 \Dc}
f) \big\|_{L^p(\mu_\chi)} \\
& \lesssim \big\|e^{-\frac 14 \Dc} f\big\|_{L^p(\mu_\chi)} \,,
\end{align*}
and thus \eqref{Besov-TL-ineq-1} holds true.\ms

In order to proceed with the main part of the 
estimates, we need to consider, for $m>\alpha/2$,
\begin{align}
W^{(m)}_t \Dc^{\gamma/2} f
& = \frac{1}{\Gamma(k-\gamma/2)} \bigg( \int_0^1 +\int_1^{+\infty} \bigg)
s^{k-\gamma/2}e^{-s\Dc} \Dc^k W^{(m)}_t f\, \frac{ds}{s} \notag \\
& =:
F_1(t,\cdot) +F_\infty(t,\cdot)  \label{F-1+F-infty}\,,
\end{align}
and observe that both $F_1,F_\infty$ are functions defined on
$(0,1)\times G$. 
\ms

{\em Step 2.} We first prove that
$$
\Dc^{\gamma/2}: B^{p,q}_{\alpha+\gamma}\to B^{p,q}_\alpha
$$
 is bounded, by showing that for any $f\in\cS$
\begin{equation}\label{Besov-ineq-2}
\int_0^1 \Big( t^{-\alpha/2} \big\| W^{(m)}_t \Dc^{\gamma/2} f
\big\|_{L^p(\mu_\chi)} \Big)^q \frac{dt}{t} 
\lesssim \int_0^1 \Big( t^{-(\alpha+\gamma)/2} \big\| W^{(m+k)}_t  f
\big\|_{L^p(\mu_\chi)} \Big)^q \frac{dt}{t} \,. 
\end{equation}
By the norm equivalence  in Theorem~\ref{teo-equiv1} and~\eqref{Besov-TL-ineq-1}, 
Step 2 will follow. \ms

We now prove \eqref{Besov-ineq-2}. Using the decomposition \eqref{F-1+F-infty}, we first estimate the
latter term.  We observe that, by the Cauchy--Schwarz inequality and the
boundedness of the $g$-function~\eqref{g-fct-est}
\begin{align*}
\| F_\infty(t,\cdot)\|_{L^p(\mu_\chi)}
& \lesssim \Big\| 
\int_1^{+\infty} s^{k-\gamma/2}| e^{-s\Dc} \Dc^k W^{(m)}_t f| \, \frac{ds}{s} 
\Big\|_{L^p(\mu_\chi)} \\
& \le \Big\|  \bigg( 
\int_1^{+\infty} |(s\Dc)^k e^{-s\Dc} (W^{(m)}_t f)|^2\, \frac{ds}{s}  \bigg)^{1/2}
\Big\|_{L^p(\mu_\chi)} \\
& \lesssim  \big\| g_k(W^{(m)}_t f)  \big\|_{L^p(\mu_\chi)} \\
& \lesssim   \big\| W^{(m)}_t f  \big\|_{L^p(\mu_\chi)}  \,.
\end{align*}
Therefore,
\begin{align}\label{F-infty-est-Besov}
\int_0^1 \Big( t^{-\alpha/2} \big\|  F_\infty(t,\cdot)
\big\|_{L^p(\mu_\chi)} \Big)^q \frac{dt}{t} 
& \lesssim  \int_0^1 \Big( t^{-\alpha/2} \big\|  W^{(m)}_t f
\big\|_{L^p(\mu_\chi)} \Big)^q \frac{dt}{t} \lesssim \|
f\|_{B^{p,q}_\alpha} \lesssim \|
f\|_{B^{p,q}_{\alpha+\gamma}}\,.  
\end{align}
Next,
\begin{align}
\| F_1(t,\cdot)\|_{L^p(\mu_\chi)}
& \lesssim \int_0^1 s^{k-\gamma/2} \big\| e^{-s\Dc} \Dc^k W^{(m)}_t
f\|_{L^p(\mu_\chi)}  \,
\frac{ds}{s} \notag\\
& = \bigg( \int_0^t +\int_t^1\bigg) t^{-k} s^{k-\gamma/2} \big\|
e^{-s\Dc} W^{(m+k)}_t f \big\|_{L^p(\mu_\chi)}  \,
\frac{ds}{s} \notag \\
& =: I(t) +I\!I(t)\,. \label{F-1-decom-Besov}
\end{align}
Now,
\begin{align*}
I(t)
& = \int_0^t  t^{-k} s^{k-\gamma/2} \big\|
e^{-s\Dc} W^{(m+k)}_t f \big\|_{L^p(\mu_\chi)}  \,
\frac{ds}{s} \lesssim t^{-k} 
\| W^{(m+k)}_t f\|_{L^p(\mu_\chi)} \int_0^t s^{k-1-\gamma/2}\, ds \\
& \approx t^{-\gamma/2} 
\| W^{(m+k)}_t f\|_{L^p(\mu_\chi)}\,,
\end{align*}
so that
\begin{align}
\int_0^1 \Big( t^{-\alpha/2} I(t) \Big)^q \frac{dt}{t} 
& \lesssim  \int_0^1 \Big( t^{-(\alpha+\gamma)/2} \big\|  W^{(m+k)}_t f
\big\|_{L^p(\mu_\chi)} \Big)^q \frac{dt}{t} \lesssim \|
f\|_{B^{p,q}_{\alpha+\gamma}} \,.
\label{I-est-Besov}
\end{align}
Next,
\begin{align*}
t^{-\alpha/2} I\!I(t)
& = \int_t^1  t^{-k-\alpha/2} s^{k-\gamma/2} \big\|
e^{-s\Dc} W^{(m+k)}_t f \big\|_{L^p(\mu_\chi)}  \,
\frac{ds}{s} \\
& = \int_0^1 {\mathbf 1}_{\{ t<s\}} \Big( \frac st \Big)^{k+\alpha/2}
s^{-(\alpha+\gamma)/2} \big\| e^{-s\Dc} W^{(m+k)}_t
f\big\|_{L^p(\mu_\chi)} \, ds \\
& = \int_0^1 {\mathbf 1}_{\{ t<s\}} \Big( \frac ts \Big)^{m-\alpha/2}
s^{-(\alpha+\gamma)/2} \big\| e^{-t\Dc} W^{(m+k)}_s
f\big\|_{L^p(\mu_\chi)} \, ds \,.
\end{align*}
Hence,
\begin{align*}
\int_0^1 \Big( t^{-\alpha/2} I(t) \Big)^q\, \frac{dt}{t} 
& =:  \int_0^1  \bigg(
\int_0^1 K(s,t) 
|g(s)| \, \frac{ds}{s} \bigg)^q \,\frac{dt}{t} \,,
\end{align*}
where 
$$
K(s,t) = {\mathbf 1}_{\{ t<s\}} \Big( \frac ts \Big)^{m-\alpha/2}
\quad
\text{and}
\quad
g(s) = s^{-(\alpha+\gamma)/2} \big\| e^{-s\Dc} W^{(m+k)}_s
f\big\|_{L^p(\mu_\chi)} \,.
$$
It is easy to check that 
$$
\int_0^1 K(s,t)\, \frac{ds}{s}\lesssim 1 \quad \text{and}\quad \int_0^1 K(s,t)\,
\frac{dt}{t}  \lesssim 1
$$
so that Schur's lemma (see \cite{Folland} e.g.) gives that
\begin{align}
\int_0^1 \Big( t^{-\alpha/2} I\!I(t) \Big)^q\, \frac{dt}{t} 
& \lesssim \int_0^1  
|g(s)|^q \, \frac{dt}{t} \,. \label{I-est-Besov}
\end{align}
Thus, putting together \eqref{F-infty-est-Besov} to \eqref{I-est-Besov}
we obtain \eqref{Besov-ineq-2}.
This completes Step 2.
\ms

{\em Step 3.}  
 We now prove that 
$$
\Dc^{\gamma/2}: F^{p,q}_{\alpha+\gamma}\to F^{p,q}_\alpha
$$
is bounded, by showing that for any $f\in\cS$
\begin{equation}\label{T-L-ineq-2}
\Big\| \bigg( \int_0^1 \Big( t^{-\alpha/2} \big| W^{(m)}_t \Dc^{\gamma/2} f
\big| \Big)^q  \,\frac{dt}{t} \bigg)^{1/q} \Big\|_{L^p(\mu_\chi)}
\lesssim \Big\| \bigg( \int_0^1 \Big( t^{-(\alpha+\gamma)/2} \big| W^{(m)}_t  f
\big| \Big)^q \,\frac{dt}{t} \bigg)^{1/q} \Big\|_{L^p(\mu_\chi)}
\,.
\end{equation}
Again, this, together with Theorem~\ref{teo-equiv1} and~\eqref{Besov-TL-ineq-1},
 will give the desired boundedness.\ms

We use decomposition \eqref{F-1+F-infty}
again.  
We notice that, since $f\in\cS$ we can switch the integration order so
that 
\begin{align*}
F_\infty(t,\cdot)
& = \frac{1}{\Gamma(k-\gamma/2)} t^m e^{-t\Dc}  \int_1^{+\infty} 
s^{k-\gamma/2}e^{-(s-1/2)\Dc} \Dc^{m+k} e^{-\frac12\Dc} f\, \frac{ds}{s} \\
& = \frac{1}{\Gamma(k-\gamma/2)} t^m e^{-t\Dc}  \int_{\frac12}^{+\infty} 
\big(s+{\textstyle{\frac12}}\big )^{k-1-\gamma/2}e^{-s\Dc} \Dc^{m+k}
e^{-\frac12\Dc} f \, 
ds \,.
\end{align*}
Hence,
 by the Cauchy--Schwarz inequality and recalling \eqref{def-g}, we have
\begin{align*}
|F_\infty(t,\cdot)|
& \lesssim  t^m e^{-t\Dc} \bigg( \int_0^{+\infty} 
\big|(s\Dc)^{m+k} e^{-s\Dc}   e^{-\frac12\Dc} f\big|^2\, \frac{ds}{s}
\bigg)^{1/2}\\
&  =
t^m e^{-t\Dc} g_{m+k} ( e^{-\frac12\Dc} f) \,.
\end{align*}
Now we use Proposition 3.6 (ii) in \cite{BPV} and the boundedness of
the $g$-function to  obtain
\begin{align}
\Big\| \bigg( \int_0^1 \Big( t^{-\alpha/2} |F_\infty(t,\cdot)| \Big)^q \,
\frac{dt}{t} \bigg)^{1/q}  \Big\|_{L^p(\mu_\chi)} 
& \lesssim 
\Big\| \bigg( \int_0^1 \Big( t^{m-\alpha/2}  e^{-t\Dc} g_{m+k} ( e^{-\frac12\Dc} f)\Big)^q \,
\frac{dt}{t} \bigg)^{1/q}  \Big\|_{L^p(\mu_\chi)}  \notag \\ 
&  \lesssim \| g_{m+k} ( e^{-\frac12\Dc} f) \|_{L^p(\mu_\chi)} \notag \\
& \lesssim \| e^{-\frac12\Dc} f \|_{L^p(\mu_\chi)}\,. \label{F-infty-est-TL}
\end{align}
Finally,
\begin{align}
\int_0^1 \Big( t^{-\alpha/2} F_1(t,\cdot)\Big)^q \, \frac{dt}{t}
& \le \int_0^1 \bigg(t^{-\alpha/2}   \Big( \int_0^t +\int_t^1\Big)  
s^{k-\gamma/2}  \big| e^{-s\Dc}  \Dc^k W^{(m)}_t f\big|\,
  \frac{ds}{s} 
\bigg)^q \,
\frac{dt}{t} \notag\\
& =: I +I\!I\,, \label{F-1-decom-TL}
\end{align}
where in this case, $I$ and $I\!I$ are functions on $G$.
Similarly to the argument in Step 1, we have
\begin{align*}
I
& = \int_0^1 \Big(t^{-\alpha/2}   \int_0^t  
s^{k-\gamma/2}  \big| e^{-s\Dc}  \Dc^k W^{(m)}_t f\big|\,
  \frac{ds}{s} 
\Big)^q \,
\frac{dt}{t} \notag \\
& = \int_0^1 \Big(  \int_0^t  
\frac{t^{m-\alpha/2} s^{k-\gamma/2} }{  (s+t)^{k+m}} \big| W^{(m+k)}_{s+t} f\big|\,
  \frac{ds}{s} 
\Big)^q \,
\frac{dt}{t}\notag \\
& = \int_0^1 \Big(  \int_t^{2t}  
\frac{t^{m-\alpha/2} (\tau-t)^{k-1-\gamma/2} }{
  \tau^{k+m-1}} 
\big| W^{(m+k)}_\tau f\big|\,
  \frac{d\tau}{\tau} 
\Big)^q \,
\frac{dt}{t}\notag \\ 
& =   \int_0^1  \Big(  \int_0^1 K(\tau,t) g(\tau)\, \frac{d\tau}{\tau} 
\Big)^q \,
\frac{dt}{t}\,,
\end{align*}
where
$$
K(\tau,t) ={\mathbf 1}_ { \{ t<\tau<2t \} } \frac{t^{m-\alpha/2} (\tau-t)^{k-1-\gamma/2} }{
  \tau^{k+m-1-(\alpha+\gamma)/2}} 
\quad\text{and}\quad 
g(\tau)=  \tau^{-(\alpha+\gamma)/2}\big|  W^{(m+k)}_\tau f\big|\,. 
$$
Since
$$
\int_0^1  K(\tau,t)\, \frac{d\tau}{\tau} 
\lesssim 1  \quad\text{and}\quad 
\int_0^1  K(\tau,t)\,  \frac{dt}{t} \lesssim 1\,,
$$
we obtain that
\begin{align}
I 
& \lesssim \int_0^1 \Big(  
 t^{-(\alpha+\gamma)/2}\big|  W^{(m+k)}_t  f\big| 
\Big)^q \,
\frac{dt}{t}
\,.  \label{I-est-TL} 
\end{align}
To estimate $I\!I$ we argue in a different fashion, using a discretization that is at
the foundation of the norm equivalence Theorem~\ref{teo-equiv2}.  Namely, by Proposition 3.6 (iv) in \cite{BPV}
\begin{align}
\| I\!I^{1/q} \|_{L^p(\mu_\chi)}  
& = \Big\|  \bigg( \int_0^1 \Big(t^{m-\alpha/2}   \int_t^1   
s^{-m-\gamma/2}  \big| e^{-t\Dc}  W^{(m+k)}_s f\big|\,
  \frac{ds}{s} 
\Big)^q \,
\frac{dt}{t} \bigg)^{1/q} \Big\|_{L^p(\mu_\chi)} \notag \\
&\le \Big\|  \bigg( \int_0^1 \Big(t^{m-\alpha/2} e^{-t\Dc}  \int_t^1   
s^{-m-\gamma/2}  \big|    W^{(m+k)}_s f\big|\,
  \frac{ds}{s} 
\Big)^q \,
\frac{dt}{t} \bigg)^{1/q} \Big\|_{L^p(\mu_\chi)} \notag \\
&\lesssim  \Big\|  \bigg( \int_0^1 \Big(t^{m-\alpha/2}   \int_t^1   
s^{-m-\gamma/2}  \big| W^{(m+k)}_s f\big|\,
  \frac{ds}{s} 
\Big)^q \,
\frac{dt}{t} \bigg)^{1/q} \Big\|_{L^p(\mu_\chi)} \notag \\
& =: \Big\|  \bigg( \int_0^1 \Big(  \int_0^1   K(s,t)
g(s)\,
  \frac{ds}{s} 
\Big)^q \,
\frac{dt}{t} \bigg)^{1/q} \Big\|_{L^p(\mu_\chi)}
\,, \label{est-post-FS-vector-valued}
\end{align}
where we have set
$$
K(s,t) ={\mathbf 1}_ { \{ t<s\} } \Big( \frac ts \Big)^{m-\alpha/2}
\quad\text{and}\quad 
g(s)=  s^{-(\alpha+\gamma)/2}\big|  W^{(m+k)}_s f\big|\,. 
$$
Again, since 
$$
\int_0^1  K(s,t)\, \frac{ds}{s} 
\lesssim 1  \quad\text{and}\quad 
\int_0^1  K(s,t)\,  \frac{dt}{t} \lesssim 1\,,
$$
we see that
\begin{align}
\| I\!I^{1/q} \|_{L^p(\mu_\chi)}  
& \lesssim  \Big\|  \bigg( \int_0^1 \Big(  
 t^{-(\alpha+\gamma)/2}\big|  W^{(m+k)}_t  f\big| 
\Big)^q \,
\frac{dt}{t} \bigg)^{1/q}  \Big\|_{L^p(\mu_\chi)} 
\,. \label{II-est-TL} 
\end{align}
Therefore, \eqref{F-1-decom-TL}, \eqref{I-est-TL} and \eqref{II-est-TL}  give 
\begin{align}
\Big\| \bigg( \int_0^1 \Big( t^{-\alpha/2} F_1(t,\cdot)\Big)^q \,
\frac{dt}{t} \bigg)^{1/q} \Big\|_{L^p(\mu_\chi)} 
& \lesssim 
\big\| I^{1/q} \big \|_{L^p(\mu_\chi)}  +\big\| I\!I^{1/q} \big
\|_{L^p(\mu_\chi)} \notag \\
& \lesssim
\Big\| \int_0^1 \Big(  
 t^{-(\alpha+\gamma)/2}\big|  W^{(m+k)}_t  f\big| 
\Big)^q \,
\frac{dt}{t}\Big\|_{L^p(\mu_\chi)}  \notag \\
& \lesssim \|
f\|_{F^{p,q}_{\alpha+\gamma}} \,. \label{F-1-est-TL}
\end{align}
This, 
together with \eqref{F-infty-est-TL}, prove \eqref{T-L-ineq-2}, 
and finally
estimates \eqref{Besov-TL-ineq-1}  and \eqref{T-L-ineq-2} complete
Step 3. 
 \ms

{\em Step 4.} We finally consider the case $c>0$. Precisely, we show
that given any $c>0$ 
if $p,q\in (1,+\infty)$, $\alpha,\gamma\ge0$, then
$$
(\Dc+cI)^{\gamma/2} : X^{p,q}_{\alpha+\gamma}\to X^{p,q}_\alpha
$$
is bounded.
\ms

We use the subordination \eqref{subord-1} (with $n\ge [\gamma/2]+1$ in
place of $k$) 
to write 
\begin{align*}
(\Dc+cI)^{\gamma/2} f
& = \frac{1}{\Gamma(n-\gamma/2)} \int_0^{+\infty} s^{n-\gamma/2} 
e^{-s(c+\Dc)} (\Dc+cI)^n f\,\frac{ds}{s} \\
& = \sum_{k=0}^n \sigma_k \int_0^{+\infty} s^{n-\gamma/2} 
e^{-s(c+\Dc)} \Dc^k f\,\frac{ds}{s} \,,
\end{align*}
for suitable positive constants $\sigma_k$.  Therefore, in order to
estimate $W^{(m)}_t (\Dc+cI)^{\gamma/2} f$, it suffices to estimate
each term, for $k=0,1,\dots,n$,
$$
\sigma_k \int_0^{+\infty} s^{n-\gamma/2} 
e^{-cs} |e^{-s\Dc} \Dc^k f|\,\frac{ds}{s} 
\,. 
$$
However, since for all
$k=0,1,\dots,n$, $s^n e^{-cs} \lesssim s^k$ for $s\in(0,+\infty)$,
such estimates follow at once from the previous Steps 2 and 3. The
proof of the lemma is complete. \ms
\qed

We are finally ready to prove the main result of this section.

\begin{thm}\label{iso-prop}
Let $\alpha\ge 0$, $\gamma\ge0$, $p,q\in(1,+\infty)$,  and let
$X^{p,q}_\alpha$ denote either space $F^{p,q}_\alpha$ or $B^{p,q}_\alpha$.   Then for $c>0$ sufficiently large,
$$
(\Dc+cI)^{\gamma/2} : X^{p,q}_{\alpha+\gamma} \to X^{p,q}_{\alpha} 
$$
is a surjective isomorphism, and its inverse is
$(\Dc+cI)^{-\gamma/2}$.  Moreover, if $c>0$ is sufficiently
large, then for all
$\alpha,\gamma\ge0$, the operators
$$
\Dc^\gamma (\Dc+cI)^{-\gamma}: X^{p,q}_{\alpha} \to X^{p,q}_{\alpha} 
$$
are bounded.
\end{thm}

We point out that for any $\gamma>0$, $c$ 
is to be chosen so that the local Riesz
transforms 
$X_J(\Dc+cI)^{-|J|/2}$ are bounded on
$L^p(\mu_\chi)$, for $1<p<\infty$ and $|J|\le [\gamma/2]+1$.
 Moreover, the operators 
$\Dc^\gamma (\Dc+cI)^{-\gamma}$ can be thought as a simplified version
of the local Riesz transforms.

\proof
{\em Step 1.} We first prove that
for $n\in\bbN$, 
$$
(\Dc+cI)^{-n} : X^{p,q}_\alpha \to X^{p,q}_{\alpha+2n} 
$$
is bounded.
\ms

Since $(\Dc+cI)^{-\beta}: L^p(\mu_\chi)\to L^p(\mu_\chi)$ is bounded,
for $p\in(1,+\infty)$ and $\beta>0$, see \cite{Komatsu} or also
\cite{BPTV}, 
we trivially
have that
\begin{equation}\label{trivial-1}
\big\| e^{-\frac12 \Dc} (\Dc+cI)^{-n} f\big\|_{L^p(\mu_\chi)} 
\lesssim \big\| e^{-\frac12 \Dc} f\big\|_{L^p(\mu_\chi)} \,.
\end{equation}
Let $m\ge [\alpha/2]+1$. In the case of Besov spaces, it suffices to apply Theorem~3.2
in~\cite{BPTV} to obtain
\begin{align*}
\big\| W^{(m+n)}_t  (\Dc+cI)^{-n} f \big\|_{L^p(\mu_\chi)} 
& = \big\| (t\Dc)^n  (\Dc+cI)^{-n} W^{(m)}_t f
\big\|_{L^p(\mu_\chi)}  \\
& \lesssim t^n \sum_{|J|\le 2n} \big\|  X_J(\Dc+cI)^{-n} W^{(m)}_t f
\big\|_{L^p(\mu_\chi)} \\
& \lesssim t^n \big\|  W^{(m)}_t f\big\|_{L^p(\mu_\chi)} \,.
\end{align*}
Hence,
\begin{align}
\bigg( \int_0^1 \Big( t^{-(n+\alpha/2)} \big\|W^{(m+n)}_t  (\Dc+cI)^{-n}
f \big\|_{L^p(\mu_\chi)} \Big)^q \, \frac{dt}{t}
\bigg)^{1/q} 
& \lesssim \bigg( \int_0^1 \Big( t^{-\alpha/2} 
\big\|  W^{(m)}_t f\big\|_{L^p(\mu_\chi)} \Big)^q \, \frac{dt}{t} 
\bigg)^{1/q} \smallskip \notag \\
& \lesssim \| f\|_{B^{p,q}_\alpha} \,. \label{besov-2}
\end{align}
Therefore, \eqref{trivial-1} and \eqref{besov-2} show that
$(\Dc+cI)^{-n}:   B^{p,q}_\alpha \to B^{p,q}_{\alpha+2n} $ 
is bounded, $p,q\in(1,+\infty)$, $\alpha\ge0$, $n\in\bbN$.\ms

Next we consider the case of the Triebel--Lizorkin spaces.  Arguing as in \eqref{subord-1} we write 
\begin{align*}
W^{(m+n)}_t (\Dc+cI)^{-n} f
& = \frac{1}{\Gamma(n)} \int_0^{+\infty} s^n e^{-cs} e^{-s\Dc}  W^{(m+n)}_t
f\, \frac{ds}{s} \\
& = \frac{1}{\Gamma(n)} \bigg( \int_0^1 + \int_1^{+\infty} \bigg) s^n e^{-cs} e^{-s\Dc}  W^{(m+n)}_t
f\, \frac{ds}{s} \\
& =:  \frac{1}{\Gamma(n)} \big(A_1(t,\cdot) +   A_\infty(t,\cdot)\big)
 \,.
\end{align*}
We begin with the latter term and observe that, since $f\in\cS$ we can switch
the integration order, so that by the Cauchy--Schwarz inequality
\begin{align*}
 |A_\infty(t,\cdot) |
& = t^{m+n}\Big| e^{-t\Dc} \int_1^{+\infty}  s^n e^{-cs} e^{-(s-1/2)\Dc}  \Dc^{m+n}
e^{-\frac12\Dc} f\, \frac{ds}{s}\Big| \\
& \lesssim t^{m+n} e^{-t\Dc}  \bigg(
\int_0^{+\infty}  \big| W^{(m+n)}_s  
e^{-\frac12\Dc} f \big|^2 \, \frac{ds}{s} \bigg)^{1/2} \\
& =  t^{m+n} e^{-t\Dc} g_{m+n}(e^{-\frac12\Dc} f)
\,,
\end{align*}
where $g_{m+n}$ is defined in \eqref{def-g}.
Therefore, by the above estimate, Proposition 3.6 (ii) in \cite{BPV}
and the boundedness of the $g$-function,
\begin{align}
\Big\| \bigg(
\int_0^1 \Big( t^{-(n+\alpha/2)} \big| A_\infty(t,\cdot) \big|
 \Big)^q \, \frac{dt}{t}
\bigg)^{1/q} \Big\|_{L^p(\mu_\chi)} 
& \lesssim  \Big\| \bigg( \int_0^1 \Big( t^{m-\alpha/2} e^{-t\Dc}
g_{m+n}(e^{-\frac12\Dc} f)
  \Big)^q \, \frac{dt}{t}
\bigg)^{1/q} \Big\|_{L^p(\mu_\chi)} 
 \notag\\
& \lesssim \big\| g_{m+n}(e^{-\frac12\Dc} f)\big\|_{L^p(\mu_\chi)}
\notag\\
& \lesssim \big\| e^{-\frac12\Dc} f \big\|_{L^p(\mu_\chi)}
\,.  
\label{da-mettere}
\end{align}

Now we turn to $A_1(t,\cdot)$, and also in this case we can switch the
integration order so that
\begin{align*}
|A_1(t,\cdot)| 
& =  \Big| t^{m+n} e^{-t\Dc} \int_0^1 s^n e^{-cs} e^{-s\Dc}  \Dc^{m+n}
f\, \frac{ds}{s} \Big| \\
& \le  t^{m+n} e^{-t\Dc} \int_0^1 s^n \big| e^{-s\Dc}  \Dc^{m+n}f
\big|   \, \frac{ds}{s} \,.
\end{align*}
Therefore, using Proposition 3.6 (ii) in \cite{BPV} we have
\begin{align}
\int_0^1 \Big( t^{-(n+\alpha/2)} |A_1(t,\cdot)| \Big)^q \, \frac{dt}{t}  
&\lesssim 
\int_0^1 \Big(   t^{m-\alpha/2} \int_0^1  s^n  |e^{-s\Dc}
\Dc^{m+n} f| \,
\frac{ds}{s} \Big)^q  \frac{dt}{t}
\notag \\
& = \int_0^1 \Big(   t^{m-\alpha/2} \bigg( \int_0^t +\int_t^1 \bigg)  s^n  |e^{-s\Dc}
\Dc^{m+n} f| \,
\frac{ds}{s} \Big)^q  \frac{dt}{t} \,. \label{intermediate-est1}
\end{align}
In order to estimate the first term on the right hand side of
\eqref{intermediate-est1} we use 
 Jensen's inequality to see that
\begin{align}
\int_0^1    t^{q(m-\alpha/2)-1} \Big( \int_0^t  s^{n-1}  |e^{-s\Dc}
\Dc^{m+n} f| \, ds \Big)^q  \, dt 
& \le \int_0^1    s^{q(m-\alpha/2)-1} s^{qn}  |e^{-s\Dc}
\Dc^{m+n} f|^q   \, ds \notag \\
 &  = \int_0^1  \Big( t^{-\alpha/2} |W^{(n+m)}_t f| \Big)^q\,
 \frac{dt}{t} \,. \label{intermediate-est2}
\end{align}
For the second term on the right hand side of
\eqref{intermediate-est1} we use Schur's test.  Precisely, arguing as
at the end of Step 2 in the proof of Lemma \ref{frac-sublap-bdd-lem},
we have
\begin{align}
\int_0^1 \Big(   t^{m-\alpha/2} \int_t^1 s^n  |e^{-s\Dc}
\Dc^{m+n} f| \,
\frac{ds}{s} \Big)^q  \frac{dt}{t}
& = \int_0^1  \Big(   \int_0^1 {\mathbf 1}_{\{t<s\}} \Big(\frac ts
\Big)^m 
  s^{-\alpha/2}\big| W^{(m+n)}_s f\big| \frac{ds}{s} \Big)^q  \, \frac{dt}{t} \notag \\
 &  \lesssim \int_0^1  \Big( t^{-\alpha/2} |W^{(n+m)}_t f| \Big)^q\,
 \frac{dt}{t} \,. \label{intermediate-est3}
\end{align}
Thus, we obtain that
\begin{align}
\Big\| \bigg(
\int_0^1 \Big( t^{-(n+\alpha/2)} \big| A_1(t,\cdot) \big|
 \Big)^q \, \frac{dt}{t}
\bigg)^{1/q} \Big\|_{L^p(\mu_\chi)} 
& \lesssim  \Big\| \bigg( \int_0^1  \Big( t^{-\alpha/2} |W^{(n+m)}_t f| \Big)^q
 \, \frac{dt}{t}
\bigg)^{1/q} \Big\|_{L^p(\mu_\chi)} \notag \\
& \lesssim \| f\|_{F^{p,q}_\alpha} \,, \label{da-mettere2}
\end{align}
as we wished to show.  Estimates \eqref{da-mettere} and \eqref{da-mettere2}
complete the  proof of Step 1. \ms

{\em Step 2.}  We now complete the proof.\ms

We first observe that
for $\gamma\ge0$,  the operator
$$
(\Dc+cI)^{-\gamma/2} : X^{p,q}_\alpha \to X^{p,q}_{\alpha+\gamma} 
$$
is bounded. Indeed, it suffices to notice that, choosing $n\in\bbN$, $n>\gamma/2$ we have
$$
(\Dc+cI)^{-\gamma/2} = (\Dc+cI)^{n-\gamma/2}
(\Dc+cI)^{-n} \,.
$$
The conclusion follows from the boundedness of the two operators on
the right hand side given by Step~1 and Lemma~\ref{frac-sublap-bdd-lem}.
\ms

It now follows that 
$$
(\Dc+cI)^{\gamma/2} : X^{p,q}_{\alpha+\gamma} \to X^{p,q}_\alpha
$$
is a surjective isomorphisms. For given any $f\in X^{p,q}_\alpha$, indeed,
we can write
$$
f = (\Dc+cI)^{\gamma/2}  (\Dc+cI)^{-\gamma/2} f\,,
$$
since $f\in L^p(\mu_\chi)$.
\ms

Finally, it is now clear that for $\alpha,\gamma\ge0$, 
$$
\Dc^{\gamma/2}(\Dc+cI)^{-\gamma/2} : X^{p,q}_\alpha \to X^{p,q}_\alpha
$$
is bounded.  The proof is now complete.
\qed

\ms

\section{Final remarks and open problems}\label{geo-inq-sec}  
\ms

In  this final section we discuss some directions for future work and
indicate some open problems.   

\ms

First of all, we stress the fact that in this work and~\cite{BPV}
we have limited ourselves to the cases $p,q\in[1,+\infty]$ and
$\alpha\ge0$. It would be interesting to investigate whether the spaces $F^{p,2}_0(\mu_\chi)$ with $p=1,+\infty$, correspond respectively 
to the local Hardy space $\mathfrak h^1(\mu_\chi)$ and its dual $\mathfrak{bmo}(\mu_\chi)$, introduced in \cite{BPTV}, in analogy to the Euclidean setting. Such spaces 
turn out to be useful in many problems,
most noticeably in the boundedness of singular integral 
operators. Moreover, Triebel--Lizorkin and Besov spaces with $0<p,q<1$ are quasi-Banach and their treatment often requires different techniques. Finally, the spaces $X^{p,q}_\alpha(\mu_\chi)$ with $\alpha<0$ should appear as
natural duals of the spaces with positive index of regularity and are
also of considerable interest.  \ms

We recall that Besov and also
Triebel--Lizorkin spaces are instrumental to applications to solvability and regularity
of solutions of nonlinear differential equations, as, for instance, in the spirit of the results 
in Section~6 in~\cite{BPTV}.  It would also be interesting to study
the homogeneous versions of Sobolev, Besov and Triebel--Lizorkin
spaces in the setting of this work.  These spaces, in particular the
homogeneous Besov spaces, appear naturally in the Strichartz estimates
for the wave equation in the Euclidean space, or Lie groups of
polynomial growth, see e.g.~\cite{GV}, \cite{BGX} and~\cite{FMV1}. \ms

Another set of natural and interesting questions concerns the 
generalization of some classical geometric
inequalities, which have already been studied in the setting of manifolds and metric spaces under suitable geometric assumptions. 
In particular, we mention the Poincar\'e inequality, see~\cite{Leoni} for the classical case
and~\cite{Jerison-Duke86} for Carnot--Carath\'eodory groups, 
the trace inequalities, see~\cite{Leoni} for the classical case and
\cite{DGN} for Carnot--Carath\'eodory groups, isoperimetric and Sobolev
inequalities \cite{Leoni} and \cite{GN}, to name just a few. We intend to study extensions of
  these classical inequalities to the case of the sub-Laplacian
  $\Dc$ on a general Lie group $G$ and of the Sobolev, Triebel--Lizorkin and Besov spaces. We point out that in \cite{RY} the authors proved 
versions of Hardy, Hardy--Sobolev, Caffarelli--Nirenberg,
Gagliardo--Nirenberg 
inequalities in the case of the Sobolev spaces $L^p_\alpha(\mu_\chi)$.

\bibliography{proc-GAHA-bib}
\bibliographystyle{abbrv}

\end{document}